\documentclass{article}
\usepackage{amsmath,amsthm}
\usepackage{amssymb}
\usepackage{yfonts,mathrsfs}
\usepackage{color}

\usepackage{latexsym}
\usepackage[english]{babel}
\usepackage[dvips]{graphicx}
\usepackage{xypic}

\begin{document}

\newcommand{\J}{\mathcal{J}}
\newcommand{\Id}{\mathrm{Id}}
\newcommand{\nei}{neighborhood}
\newcommand{\Kil}{\mathrm{Kil}}
\newcommand{\g}{{\mathbf g}}
\newcommand{\prob}{{(\bf P)}}
\newcommand{\virg}[1]{\textquotedblleft#1\textquotedblright}

\newcommand{\matSU}{\Pt{\ba{cc}
\alpha&\beta\\
-\bar \beta& \bar \alpha\ea}}

\newcommand{\funz}[5]{#1 : \begin{tabular}{ccl}
 #2 &$\rightarrow$& #3, \\
 #4 & $\mapsto$& #5   \end{tabular}}

\newcommand{\mfunz}[5]{\funz{$#1$}{$#2$}{$#3$}{$#4$}{$#5$}}
\newcommand{\mmfunz}[5]{\begin{center}
\funz{$\displaystyle{#1}$}{$\displaystyle{#2}$}{$\displaystyle{#3}$}{$\displaystyle{#4}$}{$\displaystyle{#5}$}
\end{center}}
\newcommand{\fz}[3]{#1:\, #2 \rightarrow #3}

\newcommand{\con}{C}
\newenvironment{graffa3}{\left\{ \begin{array}{lcc}}{\end{array} \right.}

\newcommand{\immagine}[3][8]{
\begin{figure}[hbpt]
\begin{center}
\includegraphics[width=#1cm]{immagini/#2}
\caption{#3}
\label{fig:#2}
\end{center}
\end{figure}}

\newcommand{\immintro}[3][8]{
\begin{figure}[hbpt]
\begin{center}
\includegraphics[width=#1cm]{immagini/#2}
\caption{#3}
\label{intro:#2}
\end{center}
\end{figure}}

\renewcommand{\span}[1]{\mathrm{span}\Pg{#1}}
\newcommand{\EXP}{\mathrm{Exp}}
\newcommand{\Exp}{\mathrm{Exp}}
\newcommand{\Tr}{\mathrm{Tr}}
\newcommand{\Mat}{\mathrm{Mat}}
\newcommand{\Pt}[1]{\left( #1 \right)}
\newcommand{\Pg}[1]{\left\{ #1 \right\}}
\newcommand{\Pq}[1]{\left[ #1 \right] }
\newcommand{\Pa}[1]{\left\langle #1 \right\rangle }
\newcommand{\hp}{hypothesis}
\renewcommand{\Re}[1]{\mathrm{Re}\Pt{#1}}
\renewcommand{\Im}[1]{\mathrm{Im}\Pt{#1}}

\newcommand{\bbibitem}{\bibitem}

\newcommand{\llabel}[1]{{\label{#1}}}

\newcommand{\ffoot}[1]{}

\renewcommand{\r}[1]{(\ref{#1})}
\newcommand{\bi}{\begin{itemize}}
\newcommand{\ei}{\end{itemize}}
\newcommand{\bd}{\begin{description}}
\newcommand{\ed}{\end{description}}
\newcommand{\be}{\begin{enumerate}}
\newcommand{\ee}{\end{enumerate}}

\renewcommand{\i}{\item}

\newcommand{\bqn}{\begin{eqnarray}}
\newcommand{\eqn}{\end{eqnarray}}
\newcommand{\eqnn}{\nonumber\end{eqnarray}}
\newcommand{\eqnl}[1]{\label{#1}\end{eqnarray}}
\newcommand{\bcas}{\begin{cases}}
\newcommand{\ecas}{\end{cases}}

\newcommand{\nn}{\nonumber\\}
\newcommand{\ba}[1]{\begin{array}{#1}}
\newcommand{\ea}{\end{array}}

\newcommand{\R}{\mathbb{R}}
\newcommand{\C}{\mathbb{C}}

\newtheorem{ml}{\bf Lemma}
\newtheorem{Theorem}{\bf Theorem}
\newtheorem{mo}{\bf \underline{{\sl Observation}}}
\newtheorem{mrem}{\bf \underline{{\sl Remark}}}
\newtheorem{mcc}{\bf Corollary}
\newtheorem{Definition}{\bf Definition}
\newtheorem{mpr}{\bf Proposition}
\newtheorem{mproperty}{\bf Property}

\newcommand{\bt}{\begin{Theorem}}
\newcommand{\et}{\end{Theorem}}
\newcommand{\bl}{\begin{ml}}
\newcommand{\el}{\end{ml}}
\newcommand{\bp}{\begin{mpr}}
\newcommand{\ep}{\end{mpr}}
\newcommand{\bc}{\begin{mcc}}
\newcommand{\ec}{\end{mcc}}
\newcommand{\bdeff}{\begin{Definition}}
\newcommand{\edeff}{\end{Definition}}
\newcommand{\brem}{\begin{mrem}\rm}
\newcommand{\erem}{\end{mrem}}

\renewcommand{\proof}{{\bf Proof. }}
\newcommand{\Lam}{\Lambda}
\newcommand{\Ga}{\Gamma}
\newcommand{\lam}{\lambda}
\newcommand{\al}{\alpha}
\newcommand{\eps}{\varepsilon}
\newcommand{\de}{\delta}
\newcommand{\om}{\omega}
\renewcommand{\th}{\theta}
\newcommand{\Z}{{\mathbb Z}}
\renewcommand{\l}{\mbox{{\footnotesize \bf L}}}

\newcommand{\eproof}{\qed} 

\renewcommand{\k}{{\mbox{\bf k}}}
\newcommand{\p}{{\mbox{\bf p}}}

\newcommand{\ga}{\gamma}

\newcommand{\Lqu}{L(4,1)}
\newcommand{\GG}{{\g}}
\newcommand{\Dir}[1]{{\mathrm{Dir}}\Pt{#1}}

\thispagestyle{empty}
\begin{center} \noindent
{\LARGE{\sl{\bf Projective Reeds-Shepp car on $S^2$ with quadratic cost}}}
\vskip 1cm
Ugo Boscain

{\footnotesize LE2i,  CNRS  UMR5158, 
Universit\'{e} de Bourgogne,
9, avenue Alain Savary - BP 47870,   
21078 Dijon CEDEX, France}

and

{\footnotesize SISSA, via Beirut 2-4 34014 Trieste, Italy - {\tt boscain@sissa.it}}\\
\vspace*{.5cm}

Francesco Rossi\\
{\footnotesize SISSA, via Beirut 2-4 34014 Trieste, Italy - {\tt rossifr@sissa.it}}

\end{center}
\vspace{.5cm} \noindent \rm
\begin{quotation}
\noindent  {\bf Abstract}
Fix two points $x,\bar{x}\in S^2$ and two directions (without orientation) $\eta,\bar\eta$ of the velocities in these points. In this paper we are interested to the problem of minimizing the cost 
$$ J[\ga]=\int_0^T \Pt{\GG_{\gamma(t)}(\dot\gamma(t),\dot\gamma(t))+ 
 K^2_{\ga(t)}\GG_{\gamma(t)}(\dot\gamma(t),\dot\gamma(t))} ~dt$$
along all smooth curves starting from $x$ with direction $\eta$ and ending in $\bar{x}$ with direction $\bar\eta$. Here $\g$ is the standard Riemannian metric on $S^2$ and $K_\ga$ is the corresponding geodesic curvature.

\noindent
The interest of this problem comes from mechanics and geometry of vision. It can be formulated as a sub-Riemannian problem on the lens space $L(4,1)$.

\noindent
We compute the global solution for this problem: an interesting feature is that some optimal geodesics present cusps. The cut locus is a stratification with non trivial topology.
\end{quotation}


\vskip 1cm
\begin{center}
PREPRINT SISSA 34/2008/M
\end{center} 

\newpage
\section{Introduction}
Fix two points $x$ and $\bar{x}$ on  a 2-D Riemannian manifold $(M,\GG)$ and two directions  $\eta\in PT_{x} M$ and $\bar{\eta}\in 
PT_{\bar{x}}M$. Here by $PT_x M$ we mean the projective tangent space at the point $x$, i.e. the  tangent space $T_x M\backslash\Pg{0}$ with the identification $v_1\sim v_2$  if there exists $\al\in\R\setminus\{ 0\}$ such that $v_1=\al v_2$. Given a vector $v\in T_x M\backslash\Pg{0}$ we call $\Dir{v}$ its direction in $PT_xM$.
We are interested in finding the path $\gamma:[0,T]\to M$ minimizing 
a  compromize among length and geodesic curvature of a curve on $M$ with fixed initial and final points $x, \bar{x}$, and fixed initial and final directions of 
the velocity $\eta, \bar{\eta}$. More precisely we consider the minimization  problem: 
 \bqn
J[\ga]=\int_0^T \left({\GG_{\gamma(t)}(\dot\gamma(t),\dot\gamma(t))+ 
\beta^2 K^2_{\ga(t)}\GG_{\gamma(t)}(\dot\gamma(t),\dot\gamma(t)) }\right)~dt~~\to~~\min,
\eqnl{eq-KOST-1}
along all smooth curves satisfying the boundary conditions 
$\gamma(0)=x$,  
$\Dir{\dot\gamma(0)}=\eta$,  
$\gamma(T)= \bar{x}$,  
$\Dir{\dot\gamma(T)}= \bar{\eta}$, as in Figure \ref{intro:ammissibile}. Here $K_{\ga(t)}$ is the geodesic curvature of 
the curve $\gamma(t)$,  $\beta$ is a fixed  constant and the final time $T$ is fixed.
\immintro[7]{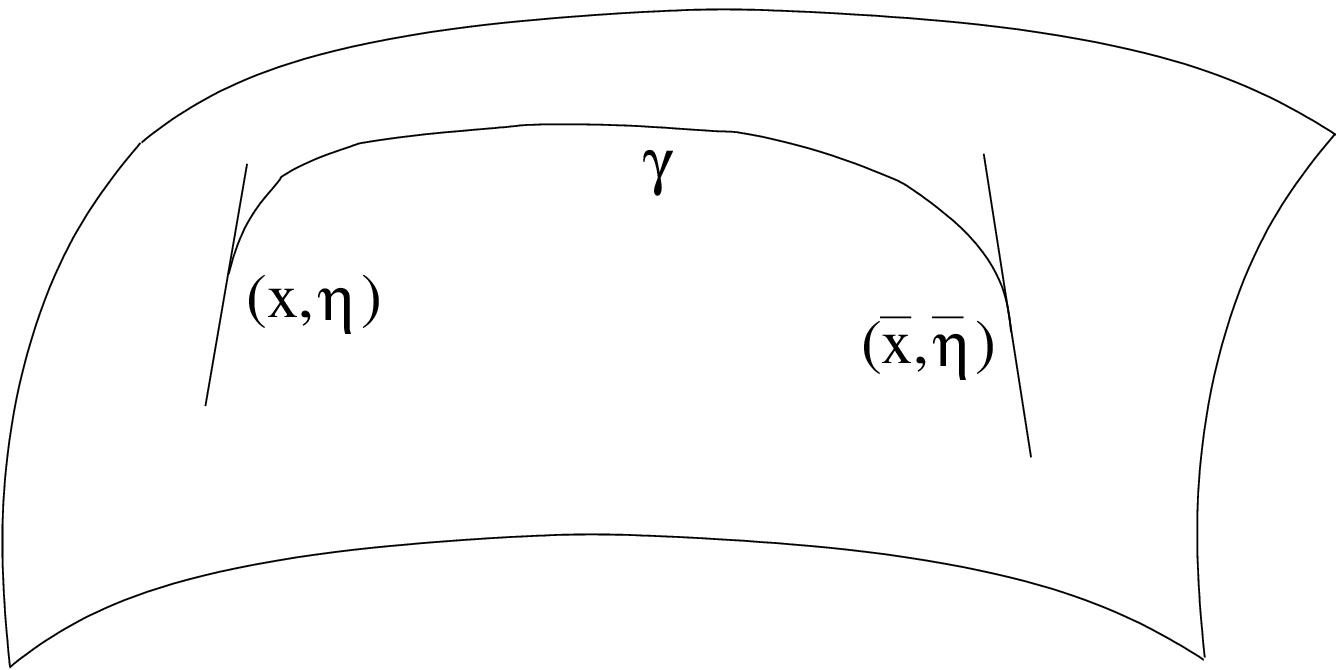}{A curve $\ga$ satisfying the boundary conditions.}

Notice that requiring $\ga$ smooth (indeed, $\ga\in C^2$ is enough) is necessary for the geodesic curvature being defined in all points. However, we will formulate the problem in a wider class to be able to apply standard existence results and the Pontryagin Maximum Principle (PMP in the following, see for instance \cite{agra-book, pontlibro}), that is a first order necessary condition for optimality.


The problem stated above is extremely difficult in general. Indeed, one can see (see Section \ref{s-mecc}) that on the projective tangent bundle  $PTM$ the minimization problem \r{eq-KOST-1} gives rise to a contact 3-D sub-Riemannian problem for which global solutions are known only for few examples. The typical difficulties one meets in such problems are the following:
\be
\i One should apply first order necessary conditions for optimality (PMP) and solve an Hamiltonian system that generically is not integrable, to find candidate optimal trajectories (called {\it geodesics}).
\i Even if all solutions to the PMP are found, one has to evaluate their optimality.
\ee

In this paper we present the solution of this problem on the sphere $S^2$, that we call the {\it projective Reeds-Shepp car with quadratic cost}. In this case  the Hamiltonian system given by the PMP is Liouville integrable and one is faced to the problem of evaluating optimality of the geodesics. This second problem is solved with the computation of the cut locus $K_x$, that is the set of points reached by more than one optimal geodesic starting from $x$. For our specific problem, indeed the cut locus coincides with the set of points in which geodesics lose optimality.

The interest of this problems comes from mechanics, and more recently from problems of geometry of vision \cite{citti-sarti,petitot,Reeds-Shepp}. Indeed, in the spirit of the model of visual architecture  due to  Petitot, Citti-Sarti and Agrachev,  the minimization problem stated above in the case $M=\R^2$ is the optimal control problem solved by the visual cortex $V1$ to reconstruct a contour that is partially hidden or corrupted in a planar image. The case $M=S^2$ can be seen as a modified version of this model, taking into account the curvature of the {\it retina} an/or the movements of the eyes.

\subsection{Examples of problems of type \r{eq-KOST-1}}
\subsubsection{The planar case}
We present the case $M=\R^2$ endowed with the standard Riemannian metric: it is being studied in \cite{igorino}. In this case the problem \r{eq-KOST-1} is equivalent to the following optimal control problem
\bqn
&&\ba{c}\left(\ba{c}
\dot x\\
\dot y\\
\dot \th
\ea
\right)=u_1\left(\ba{c}
\cos(\th)\\
\sin(\th)\\
0
\ea
\right)+u_2\left(\ba{c}
0\\
0\\
1
\ea
\right),~~~~~~~~~~
\int_0^T (u_1^2+\beta^2 u_2^2) dt~\to~\min,\\\\
(x(0),y(0),\th(0))=(x_0,y_0,\th_0),~~~~~~(x(T),y(T),\th(T))=(x_1,y_1,\th_1).\ea
\eqnn 
The dynamics coincide with the Reeds-Shepp car (see \cite{Reeds-Shepp}), except for the fact that here $\theta\in\R/\pi$. For this reason we call the problem \r{eq-KOST-1} the {\bf projective Reeds-Shepp car} on a Riemannian manifold, with quadratic cost.

\subsubsection{The case $M=S^2$}
The case $M=S^2$, endowed with the standard Riemannian metric, can be seen as the compactified version of the problem on the plane. It is the optimal control problem
\bqn
\ba{c}
\left(\ba{c}
\dot a\\
\dot b\\
\dot \xi
\ea
\right)=u_1\left(\ba{c}
\cos(\xi)\\
\displaystyle\frac{\sin(\xi)}{\sin(a)}\\
-\cot(a)\sin(\xi) 
\ea
\right)+u_2\left(\ba{c}
0\\
0\\
1
\ea
\right),~~~~~~
\int_0^T (u_1^2+\beta^2 u_2^2) dt~\to~\min,\\\\
(a(0),b(0),\xi(0))=(a_0,b_0,\xi_0),~~~~~~~(a(T),b(T),\xi(T))=(a_1,b_1,\xi_1),\ea
\eqnn 
where $(a,b)$ are standard spherical coordinates on the sphere (see \r{sfer}) and $\xi=\Dir{\dot{\ga}}.$ Also in this case the Hamiltonian system given by the Pontryagin Maximum Principle is Liouville integrable and one is faced to the problem of evaluating optimality of the geodesics.

This problem has special features. Indeed, $PTS^2$ is the 3-D manifold called {\it lens space} $\Lqu$, that can be seen as a suitable quotient of $SU(2)$ by a discrete group. Moreover, when $\beta=1$  this problem is the projection of a left-invariant sub-Riemannian problem on $SU(2)$, called $\k\oplus\p$ problem. For $SU(2)$ an explicit expression for geodesics is given and we computed the cut locus and the Carnot-Caratheodory distance in \cite{nostro-gruppi} (results are recalled in Section \ref{ss-SU2}). As a consequence, one can easily find geodesics for the problem on $PTS^2$ as projections of the geodesics for $SU(2)$. An interesting feature is that the projections of these geodesics on $S^2$ present cusps in some cases, see Figure \ref{intro:cusp}. Observe that in cusp points the tangent direction is well defined.
\immintro[7]{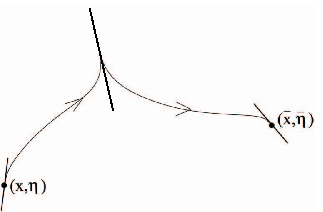}{A geodesic with a cusp.}

\newcommand{\relB}{\approx}
$\Lqu$ can be described topologically as follows:
consider a full 3-D ball $B:=\Pg{(x_1,x_2,x_3)\in\R^3\ |\ \sum_{i=1}^3 x_i\leq 1}$ and define the equivalence relation $\relB$ on it as follows: $\relB$ is reflexive; the points $(x_1^+,x_2^+,x_3^+)\in\partial B^+=\partial B\cap\Pg{x_3\geq 0}$ and $(x_1^-,x_2^-,x_3^-)\in\partial B^-=\partial B\cap\Pg{x_3\leq 0}$ are identified when $$x_3^+=-x_3^-\mbox{ and }\Pt{\begin{array}{c}
x_1^+\\
x_2^+
\end{array}}=
\Pt{\begin{array}{cc}
0 & -1\\
1 & 0
\end{array}}
\Pt{\begin{array}{c}
x_1^-\\
x_2^-\\
\end{array}},$$
as in Figure \ref{intro:L41def}. The manifold $B/_\relB$ is topologically equivalent to $\Lqu$.
\immintro[7]{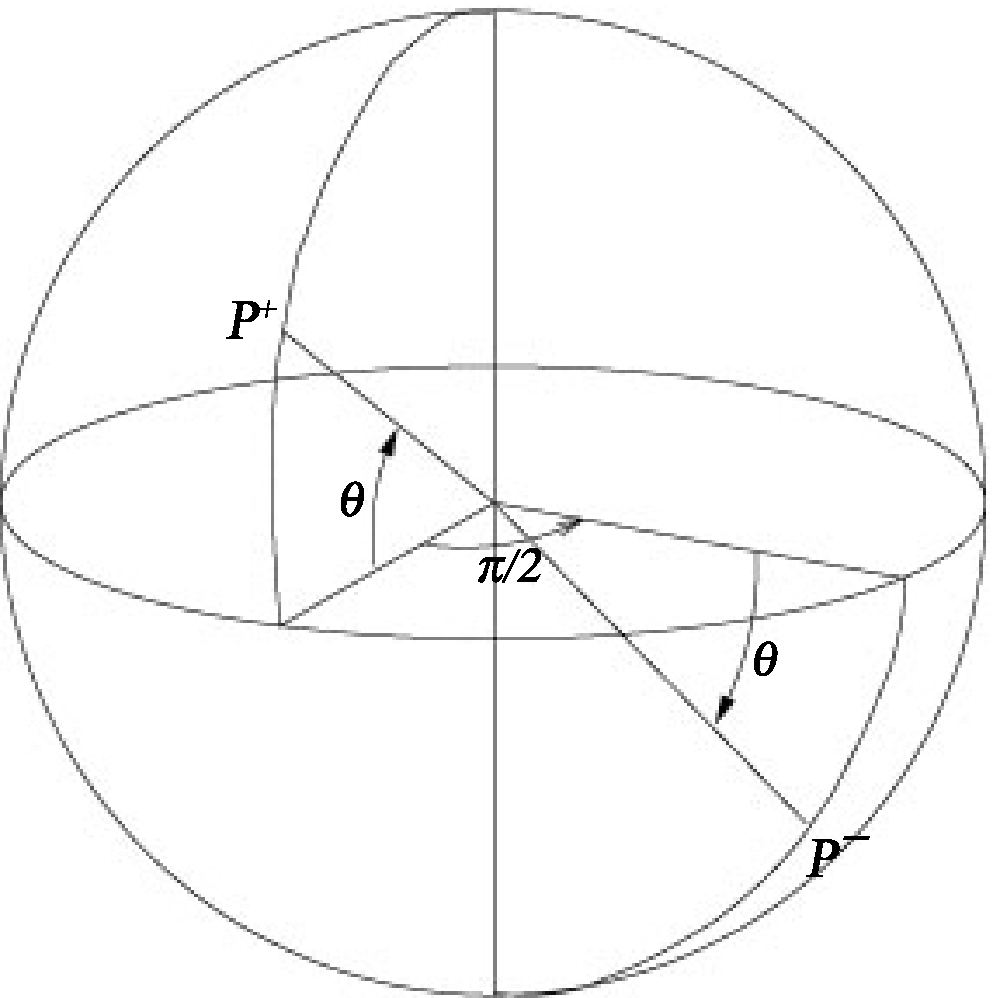}{Identification rule on $\partial B$.}

Observe that there exist different identification rules on $SU(2)$ such that the quotient manifold is topologically $\Lqu$ and the induced sub-Riemannian structure is well defined: indeed, in \cite[Section 4]{nostro-gruppi} we already endowed $\Lqu$ with a sub-Riemannian structure that is different from the one we present in this paper. The structure chosen in this paper is the unique for which the minimal length problem on $\Lqu$ coincide with the problem \r{eq-KOST-1}. As a consequence, also the cut loci are different in the two cases. See Remark \ref{r-Lqudiversi} for further details.

The main result of this paper is the computation of the cut locus $K_{[\Id]}$, presented in Theorem \ref{t-L41cut}.  It is a stratification and it has nontrivial topology. Indeed, it is the union of two sets  $K_{[\Id]}=K^{loc}\cup K^{sym}$, described as follows:
$K^{loc}$ is topologically equivalent to $S^1$ without the starting point $[\Id]$. $K^{sym}$ is given by the union of two 2-D manifolds glued along their boundaries $$K^{sym}=K^{a}\cup K^{b}$$ where
\bi
\i The manifold $K^{a}$ is contained in the set $~~~\Pg{[g]\ |\ g=\Pt{\ba{cc}\al&\beta\\-\overline{\beta}&\overline{\al}	\ea}\mbox{~with~}\Re{\beta}=0}~~~$  and it is homeomorphic to a 2-D disc with only two points identified on the boundary. Observe that the boundary is thus homeomorphic to two circles $S^1$ glued in a point (we call them $\mathscr{S}^1$ and $s^1$). A part of $K^{loc}$ (namely, the farthest part from the starting point) is contained in $K^{a}$.
\i The manifold $K^{b}$ is topologically equivalent to a 2-D disc whose boundary $S^1$ is glued to the boundary of $K^{a}$ as presented in figure \ref{intro:cutscheme}.
\ei
A representation of the topology of the cut locus is given in figure \ref{intro:cutscheme}.
\immintro[11]{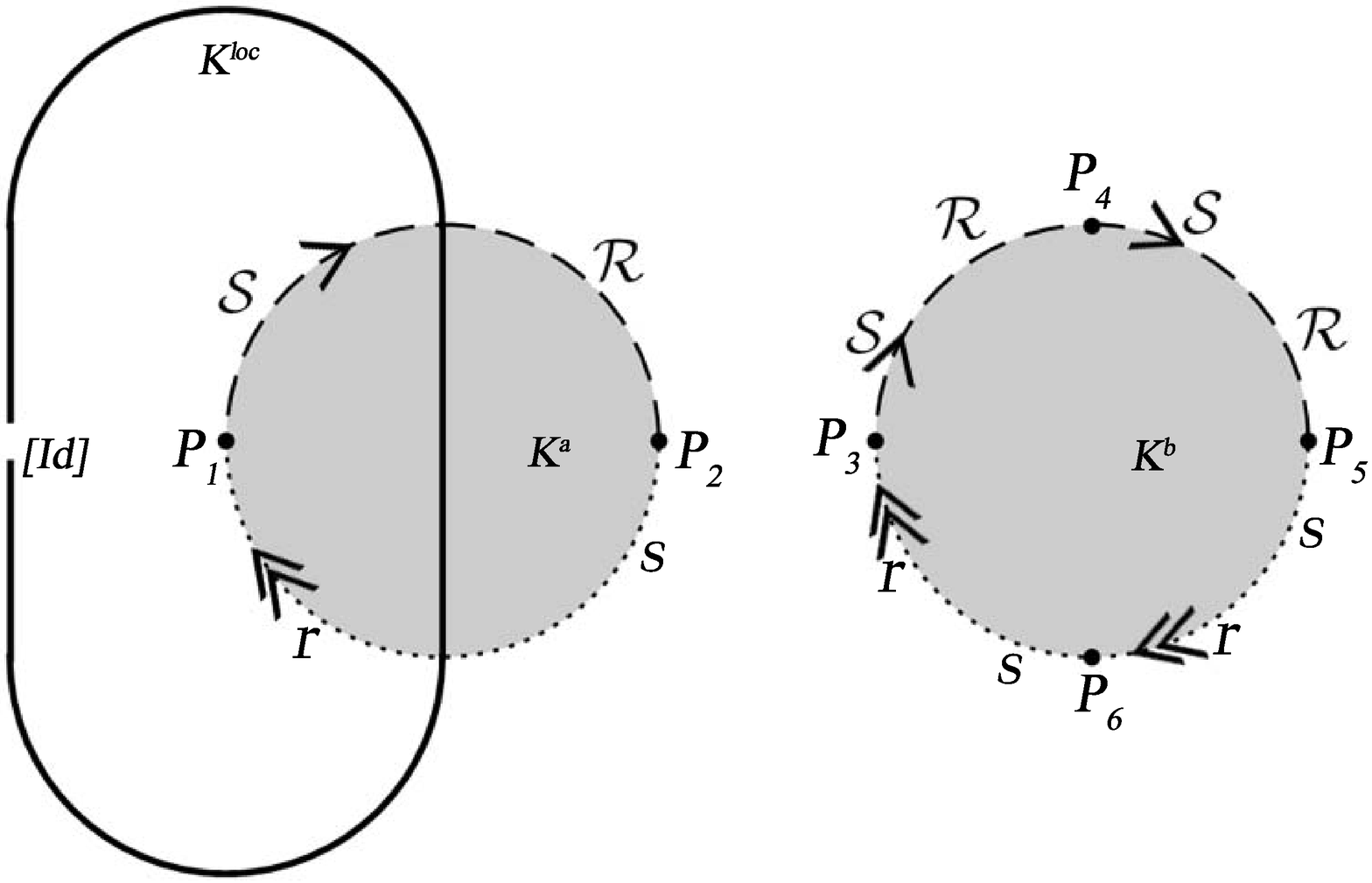}{A scheme of the topology of the cut locus: points $P_i$ are identified, boundary of disks are identified according to arrows.}

\begin{figure}[hbpt]
\begin{center}

\newcommand{\SSS}{{\cal S}}
\newcommand{\RRR}{{\cal R}}
\newcommand{\sss}{{s}}
\newcommand{\rrr}{{r}}
\scalebox{.6}{ \input{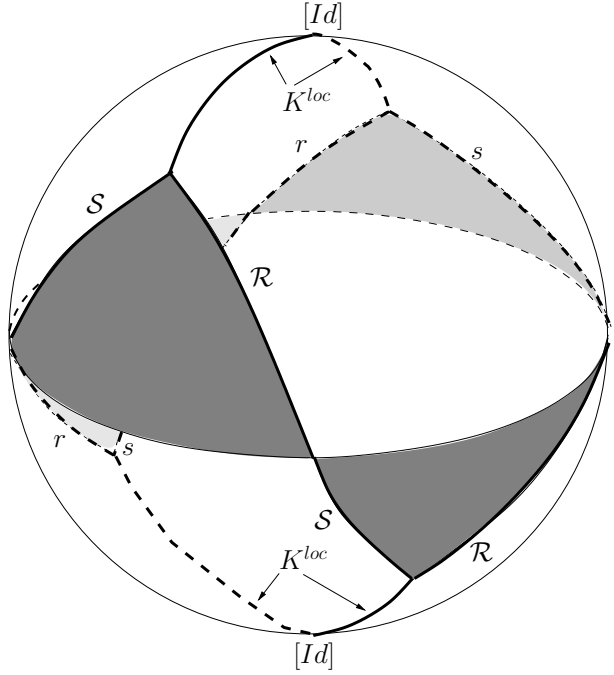}}
\caption{Pitcure of the cut locus: strata $K^{loc}$ and $K^a$.}
\label{intro:L41cut}
\end{center}
\end{figure}

We give a picture of the cut locus in Figure \ref{intro:L41cut}. Here:
\bi
\i The 1-D stratum $K^{loc}$ is represented by two semi-diameters (one is the left-right line on the northern hemisphere, the other is the front-rear line on the southern hemisphere) without the poles (representing $[Id]$): due to identification rule $\relB$, the two lines are identified.
\i The 2-D stratum $K^a$ is represented by four \virg{triangles} on the boundary of the sphere. The two dark gray triangles (on the front) are identified . Similarly, the two light gray triangles (on the rear) are identified.
\i The 2-D stratum $K^b$ is not represented for a better comprehension: it can be seen as a surface inside the sphere whose boundary is the closed line given by the concatenation of segments $\mathcal{SRSR}srsr$.
\ei

The paper is organized as follows: in Section \ref{s-mecc} we define the mechanical problem on an orientable Riemannian manifold $M$ and write the optimal control problem in the case $M=S^2$. In Section \ref{s-L41subr} we define $\Lqu$ and a sub-Riemannian structure on it such that the minimal length problem on $\Lqu$ coincide with the mechanical problem stated above. Finally, in Section \ref{s-solution} we give the solution of the problem.

\section{The problem of minimizing length and curvature on a 2-D Riemannian manifold}
\llabel{s-mecc}

In this section when we speak of {\it directions} we always mean directions without orientations.

Let $M$ be a 2-D manifold smooth manifold and $PTM$ its {\it bundle of 
directions}, i.e., 
\bqn
PTM&:=&\bigcup_{x\in M}\{ PT_xM\}, \mbox{ where }\nn\nn
PT_xM&:=&\Pg{(T_xM\setminus\{0\})/\sim \mbox{ where }v_1\sim 
v_2\mbox{ if and only if }\exists ~\al\in \R\setminus\{0\}~|~ v_1=\al v_2} .
\eqnn 

Assume now that $M$ is orientable and Riemannian, with Riemannian metric $\g$ and let  $(x_1,x_2)$ be a {\it local orthogonal chart} on an open set $U\subset M$, i.e. a chart such that,
$$\GG_{(x_1,x_2)}(\partial_{x_1},\partial_{x_2})\equiv 0 \mbox{~on~}U.$$ In the following the symbol $\|\cdot\|$ indicates the norm with respect to  
$\GG_{(x_1,x_2)}$.

In this case there is a canonical way to construct a chart 
$(x_1,x_2,\xi)$, with $\xi\in \R/\pi$ on $PTU$, as follows:
write
$$v=v_1\frac{\partial_{x_1}}{\|{\partial_{x_1}}\|}
+v_2
\frac{\partial_{x_2}}{\|{\partial_{x_2}}\|}
\in T_xU,~~x\in U, 
$$
and define $\xi$ as the angle 
between $v$ and $\partial_{x_1}/{\|{\partial_{x_1}}\|}$, i.e. $\xi:=\arctan\Pt{\frac{v_2}{v_1}}.$

In the following we call the chart $(x_1,x_2,\xi)$ a {\it 
lifted-orthogonal chart}.
Notice that  
${\|{\partial_{x_1}}\|}\equiv 1$ and ${\|{\partial_{x_2}}\|}\equiv 1$ on 
$U$, implies $\GG$ flat on $U$. \\

Given a smooth curve $\ga(.):[0,T]\to M$, its {\it lift} 
in $PTM$ is the curve $(\ga(.),\xi(.))$, where $\xi(t)\in PT_{\ga(t)}M$ is 
the direction of $\dot \ga(t)$. Writing 
$$
\dot \ga(t)=
\dot\ga_1(t){\partial_{x_1}}
+\dot\ga_2(t){\partial_{x_2}}=
(\dot\ga_1(t)
{\|{\partial_{x_1}}\|}){\partial_{x_1}}/{\|{\partial_{x_1}}\|}
+(\dot 
\ga_2(t){\|{\partial_{x_2}}\|}){\partial_{x_2}}/{\|{\partial_{x_2}}\|},
$$
then 
\bqn 
\xi(t)=\arctan\left(\frac{\dot\ga_2(t) 
\|{\partial_{x_2}}\|}{\dot \ga_1(t)\|{\partial_{x_1}}\|}\right).
\eqnl{eq-angolo}
Given a smooth curve $\Gamma(.):[0,T]\to PTM$, we say that it is {\it 
horizontal} if it is the lift of a curve on $M$. 
If $\Gamma(.)=(\ga_1(.),\ga_2(.),\xi(.))$ are its components in  
a lifted orthogonal chart, then $\Gamma(.)$ is horizontal if and only if 
condition \r{eq-angolo} holds.\\

The requirement that a smooth curve $\Gamma(.)$ in $PTM$ is horizontal, is 
equivalent to require that its velocity belongs to a 2-D distribution (see Appendix A). Let us describe this fact in detail: in a lifted orthogonal chart we have 
$\dot\Gamma(t)=\dot \ga_1(t)\partial_{x_1}+\dot 
\ga_2(t)\partial_{x_2}+\dot 
\xi(t)\partial_{\xi}$ and condition \r{eq-angolo} is equivalent to
\bqn
\dot\Gamma(t)\in \Delta(\Gamma(t)), \mbox{ where }
\Delta(x_1,x_2,\xi):=\span{F_1,F_2}, 
F_1=\left(\ba{c} \displaystyle \frac{\cos(\xi)}{\|{\partial_{x_1}}\|} \\
\displaystyle  
\frac{\sin(\xi)}{\|{\partial_{x_2}}\|}\\\Xi(x_1,x_2,\xi)  \ea\right),~~~ 
F_2=\left(\ba{c}  0\\0\\1  \ea\right). 
\eqnn
The function $\Xi(x_1,x_2,\xi)$ is an arbitrary smooth function. We choose $\Xi$ such that the vector field $F_1$ on PTM  is the {\it geodesic spray}, i.e. the infinitesimal generator of the geodesic flow (see \cite{mario,spivak}).

Notice that
$\span{F_1(x_1,x_2,\xi),F_2(x_1,x_2,\xi),[F_1,F_2](x_1,x_2,\xi)}=T_{(x_1,x_2,\xi)}(PTM)$, hence $(PTM,\Delta)$ is a 3-D contact distribution.

\subsection{A sub-Riemannian problem on $PTM$}

There is a natural sub-Riemannian problem associated with $(PTM,\Delta)$.  It comes from the problem of minimizing a compromise among length and geodesic curvature of a curve on $M$ with fixed initial and final points, and fixed initial and final directions of the velocities. More precisely we consider the following problem: 
\bd
\i[\prob]  Fix $x, \bar x\in M$, $\eta\in PT_xM$ and $\bar \eta\in 
PT_{\bar 
x}M$. 
Minimize 
\bqn
\J[\Ga]=\int_0^T \sqrt{\GG_{\gamma(t)}(\dot\gamma(t),\dot\gamma(t))+ 
\beta^2 K^2_{\ga(t)}\GG_{\gamma(t)}(\dot\gamma(t),\dot\gamma(t))} ~dt
\eqnl{eq-KOST}
along all horizontal smooth curves 
$\Gamma(.)=(\gamma(.),\xi(.)):[0,T]\to PTM$ satisfying the 
boundary conditions 
$\Gamma(0)=(x,\eta)$,  $\Gamma(T)=(\bar x,\bar \eta)$, i.e. 
$\gamma(0)=x$,  
$\Dir{\dot\gamma(0)}=\eta$,  
$\gamma(T)=\bar x$,  
$\Dir{\dot\gamma(T)}=\bar \eta$. Here $K_{\ga(t)}$ is the geodesic curvature of 
the curve $\gamma(t)$, and $\beta$ is a fixed nonvanishing constant.
\ed

Problem \prob\ is illustrated in Figure \ref{intro:ammissibile}. Its equivalence with the problem \r{eq-KOST-1} (in which the square root is absent) is eplained below. The cost \r{eq-KOST} is the most natural cost associated with $(PTM,\Delta)$. Indeed it is invariant by change of coordinates and by reparameterization of the curve $\ga$. Moreover, as it will be clear from the example below, it can be interpreted as a length in $PTM$.

The term $\g_{\gamma(t)}(\dot\gamma(t),\dot\gamma(t))$ 
corresponds to the Riemannian length on $M$, while the term   
$\beta^2 K^2_{\ga(t)}\GG_{\gamma(t)}(\dot\gamma(t),\dot\gamma(t))$ corresponds to the geodesic curvature. The presence of 
$\GG_{\gamma(t)}(\dot\gamma(t),\dot\gamma(t))$ in this second term guarantees the invariance by reparametrization.
The constant $\beta$ fixes the relative weight. It may happen that there is a natural choice for $\beta$ or even that $\beta$ is not relevant (see an example on $\R^2$ below), depending on the specific problem.

One can check  that $\GG_{\gamma(t)}(\dot\gamma(t),\dot\gamma(t))+
\beta^2 K^2_{\ga(t)}\GG_{\gamma(t)}(\dot\gamma(t),\dot\gamma(t))$ defines a positive definite quadratic form on $\Delta$, i.e.  it endows  $(PTM,\Delta)$ with a sub-Riemannian structure (see Appendix \ref{a-SubRiem}). Indeed, if $x=(x_1,x_2,\xi)$ is a lifted orthogonal chart and $\Xi$ is such that $F_1$ is the geodesic spray, we have that problem \prob\ becomes the optimal control problem:
\bqn
\dot x=u_1F_1+u_2F_2,~~
F_1=\left(\ba{c} \displaystyle \frac{\cos(\xi)}{\|{\partial_{x_1}}\|} \\
\displaystyle
\frac{\sin(\xi)}{\|{\partial_{x_2}}\|}\\\Xi(x_1,x_2,\xi)  \ea\right),~~~
F_2=\left(\ba{c}  0\\0\\1  \ea\right),~~~
\J[\Ga]=\int_0^T\sqrt{u_1^2+\beta^2u_2^2}~dt\to\min
\eqnn
that is the minimal length problem in the sub-Riemannian manifold $PTM$.

Consider now the problem of minimization of the {\it energy functional} $E[\Ga]:=\int_0^T\ {u_1^2+\beta^2u_2^2}~dt$ with fixed final time $T$ and fixed starting and ending points. If $\Ga$ is a minimizer of $E$, then it is a minimizer of $\J$ and its {\it lifted velocity} $v(t):=\sqrt{u_1^2+\beta^2u_2^2}$ is constant. On the other side, a minimizer $\Ga$ of $\J$ parametrized with constant lifted velocity $v$ is also a minimizer for $E$ with $T=\J[\Ga]/v$.
For details, see Appendix \ref{a-SubRiem}.

Finding a complete solution for the problem \prob, i.e. finding the optimal trajectory connecting each initial and final condition, is a problem of optimal synthesis in dimension 3 that is extremely difficult in general. The case $M=\R^2$ with the standard euclidean structure is being studied in \cite{igorino}. The complete solution in the case $M=S^2$ with the standard Riemannian metric is given as the solution of the sub-Riemannian problem on $\Lqu$ presented in Section \ref{s-L41subr}.

\subsubsection{The planar case}
Consider $M=\R^2$ with the standard euclidean structure. In this case $PT\R^2\simeq \R^2\times S^1$ and a lifted orthogonal chart is $(x_1,x_2,\xi)$, where $(x_1,x_2)$ are euclidean coordinates on $\R^2$. In this case, since $\|{\partial_{x_1}}\|=\|{\partial_{x_2}}\|=1$, writing $\ga(t)=(x_1(t),x_2(t))$, we have 
$$
\GG(\dot\ga,\dot\ga)=\dot x_1^2+\dot x_2^2,\mbox{~~~and~~~} 
K_{\ga}=\frac{- \dot x_2\,\ddot x_1   + 
\dot x_1\,\ddot x_2}{{\left(\dot x_1^2 + \dot x_2^2 
\right)}^{\frac{3}{2}}}=\frac{\dot\xi}{\sqrt{\dot x_1^2+\dot 
x_2^2}},\mbox{ where~~~}\xi(t)=\arctan\Pt{\frac{\dot x_2(t)}{\dot x_1(t)}}.
$$
Hence $\J[\Ga]=\int_0^T\sqrt{\dot x_1^2+\dot x_2^2+\beta^2 \dot \xi^2}~dt$. 
We choose  $\Xi(x_1,x_2,\xi)\equiv 0$, then $F_1$ is the geodesic spray and 
we have 
\bqn
&&\Delta(x_1,x_2,\xi):=\span{F_1,F_2}, \mbox{~~where~~}
F_1=\left(\ba{c} 
\cos(\xi) \\ 
\sin(\xi)\\
0\ea\right),~~~ 
F_2=\left(\ba{c}  0\\0\\1  \ea\right),~~~\J[\Ga]=\int_0^T\sqrt{u_1^2+\beta^2u_2^2}~dt.
\eqnn
The  associated optimal control problem  is 
$$
\left\{\ba{l}
\dot x_1=u_1\cos(\xi)\\
\dot x_2=u_1\sin(\xi)\\
\dot \xi=u_2,
\ea
\right.~~~u_1,u_2\in \R,~~~
\int_0^T\sqrt{u_1^2+\beta^2u_2^2}~dt\to\min.
$$
This problem can be seen as a left-invariant sub-Riemannian problem on the group $SE(2)/\sim$, the group of rototranslations of the plane
\bqn
SE(2):=\Pg{\left(
\ba{ccc}
\cos(\al)&-\sin(\al)&x_1\\
\sin(\al)&\cos(\al)&x_2\\
0&0&1
\ea
\right)\ |\ \al\in\R/2\pi,\ x_i\in\R},\eqnn
endowed with the identification rule $(\al,x_1,x_2)\sim (\al+\pi,x_1,x_2)$. In this case the constant $\beta$ can be set to 1, since  it becomes 
irrelevant after the transformation $(x_1,x_2)\to (\beta x_1,\beta x_2)$. 
This problem is being studied in \cite{igorino}.

\subsection{The case $M=S^2$}
Let $S^2=\{(y_1,y_2,y_3)\in\R^3~|~y_1^2+y_2^2+y_3^2=1\}$ and let $(a,b)$ be the 
spherical coordinates (see Figure \ref{fig:f-coord-sf}):
\bqn
\left\{
\ba{l}
y_1=\sin(a)\cos(b)\\
y_2=\sin(a)\sin(b)\\
y_3=\cos(a),
\ea
\right.~~~a\in[0,\pi],~~b\in[0,2 \pi).
\eqnl{sfer}
\immagine[6]{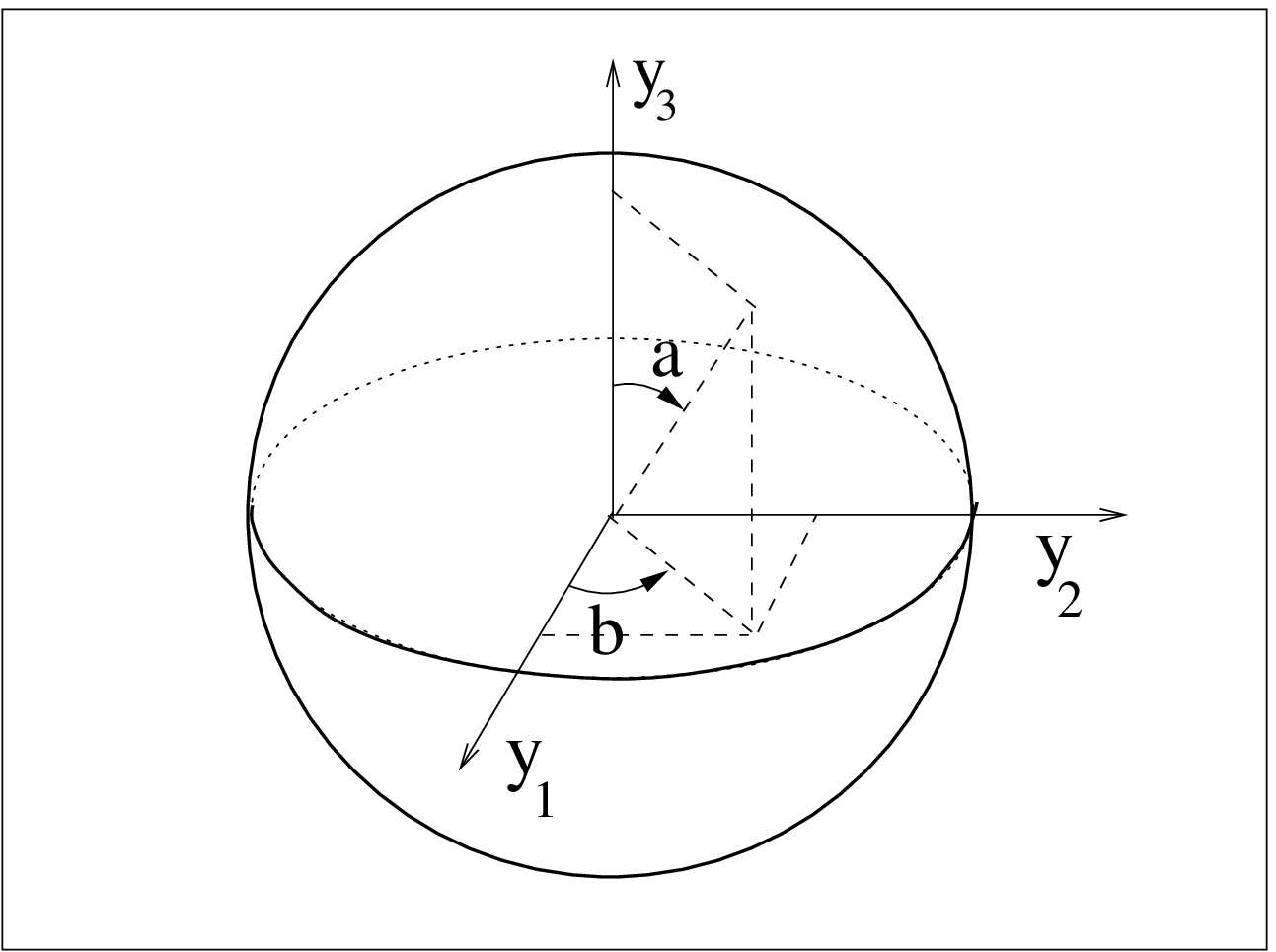}{Spherical coordinates on $S^2$.}

In these coordinates the standard Riemannian metric on $S^2$ (induced by 
its embedding in $\R^3$) has the form
\ffoot{
$$
\left\{
\ba{l}
\partial_a= 
\cos(a)\cos(b) 
\partial_{y_1}+
\cos(a)\sin(b) 
\partial_{y_2}-\sin(a)\partial_{y_3}\\
\partial_b=
-\sin(a)\sin(b) 
\partial_{y_1}+
\sin(a)\cos(b) 
\partial_{y_2}\\
\ea
\right.
$$
}
$$
\GG{(a,b)}=\left(\ba{cc}
1&0\\
0&\sin^2(a)   
\ea\right)_.
$$

The geodesic curvature of a curve $\ga(t)=(a(t),b(t))$ is:
\bqn
K_\ga&=&\frac{{\sqrt{{\sin (a)}^2}}\,\left( 2\,\cot (a)\,{\dot a}^2\,\dot 
b 
+ 
      \cos (a)\,\sin (a)\,{\dot b}^3 - \dot b\,\ddot a + \dot a\,\ddot b 
\right) }{{\left(
       {\dot a}^2 + {\sin (a)}^2\,{\dot b}^2 \right) 
}^{\frac{3}{2}}}=\frac{\cos (a)\,\dot b}{{\sqrt{\dot a^2 + {\sin 
(a)}^2\,\dot b^2}}}
+\frac{\dot \xi}{{\sqrt{\dot a^2 + {\sin (a)}^2\,\dot b^2}}}
\eqnn
with $\xi=\arctan\Pt{\frac{\dot{b}\sin(a)}{\dot{a}}}$.

Hence $\J[\Ga]=\int_0^T \sqrt{\dot a^2 + \sin^2(a)\,\dot b^2  + 
\beta^2(\cos(a)\dot b+\dot\xi)^2}$. In this case $F_1$ is the geodesic 
spray choosing $\Xi(a,b,\xi)=-\cot(a)\sin(\xi)$ for which we have 
\bqn
\Delta(a,b,\xi):=\span{F_1,F_2}, \mbox{~~where~~}
F_1=\left(\ba{c} {\cos(\xi)} \\
\displaystyle  
\frac{\sin(\xi)}{\sin(a)}\\\   -\cot(a)\sin(\xi)  
\ea\right),~~~ 
F_2=\left(\ba{c}  0\\0\\1  
\ea\right),~~~\J[\Ga]=\int_0^T\sqrt{u_1^2+\beta^2u_2^2}~dt.
\eqnn 
The  associated optimal control problem  is 
\bqn
\left\{\ba{l}
\dot a=u_1\cos(\xi)\\
\displaystyle \dot b=u_1\frac{\sin(\xi)}{\sin(a)}\\
\dot \xi=-u_1\cot(a)\sin(\xi)   +u_2,
\ea
\right.~~~u_1,u_2\in \R,~~~\int_0^T\sqrt{u_1^2+\beta^2u_2^2}~dt\to\min.
\eqnl{OCsfera}

\brem\llabel{rem-PTS2fixedpoint}
Observe that this optimal control problem admits trajectories for which $\dot{a}=\dot{b}=0$. These trajectories in $PTS^2$ represent the possibility of changing the direction of the tangent vector on the sphere $S^2$ on a fixed point; mechanically, it represents the possibility of rotation on itself.
\erem

\subsubsection{Existence of minimizers}
To discuss the problem of existence of minimizers for the problem on the sphere $S^2$, let us come back to the original problem of minimizing the cost 
\r{eq-KOST-1}, where $M=S^2$ with the standard Riemannian metric,  along all smooth curves $\gamma:[0,T]\to S^2$. Here the class of smooth curves has been chosen to give a meaning to the geodesic curvature in the whole interval $[0,T]$ (indeed $C^2([0,T])$ was enough).

Unfortunately the class of smooth (or $\con^2$) curves is too small to apply standard existence theorems for minimizers and first order necessary conditions for optimality (the PMP). We have to deal with a larger class of function.

This class arises naturally after the lift in $PTS^2$, where the minimization problem takes the form \r{OCsfera}. Here it is natural to look for minimizers $(a(.),b(.),\xi(.))$ that are absolutely continuous and corresponding to controls $u_1(.),u_2(.)$ belonging to $L^2([0,T])$. In this class, standard existence theorems can be applied to guarantee existence of minimizers (see for instance \cite[Theorem 5.1]{yuri}). 

However the PMP works in the smaller class of $L^\infty$ controls (corresponding to Lipschitz curves $(a(.),b(.),\xi(.))$); hence, before applying PMP, one needs to prove that minimizers belong to this smaller class.

Thanks to a theorem of Sarychev and Torres (see \cite[Theorem 5.2]{yuri}), one can prove that for the problem \r{OCsfera} all minimizers satisfy the conditions given by the PMP and that normal ones correspond to controls belonging to  $L^\infty([0,T])$. Since in our case there are no abnormal extremals (the problem \r{OCsfera}  is 3-D contact), it follows that all optimal controls are indeed $L^\infty([0,T])$.

Finally, all trajectories satisfying the PMP are solutions of an Hamiltonian system, that in our case is analytic. Hence all minimizers are analytic and therefore smooth. 

We conclude that it is equivalent to solve the original problem in the class of smooth curves or the lifted problem in the class of $L^2$ controls. A similar treatment has been presented with all details in \cite{yuri} for a similar problem but in the presence of a  drift.

\section{A sub-Riemannian problem on $L(4,1)$}
\llabel{s-L41subr}
In this section we define a $\k\oplus\p$ sub-Riemannian structure on $SU(2)$ and compute its cut locus and the sub-Riemannian distance. For details on Sub-Riemannian geometry and $\k\oplus\p$ manifolds, see Appendix \ref{a-SubRiem}. We then define a sub-Riemannian structure on $\Lqu$ induced by the one on $SU(2)$ and prove that the sub-Riemannian problem on $\Lqu$ coincide with the optimal control problem \r{OCsfera} on $PTS^2$.

\subsection{The $\k\oplus\p$ problem on $SU(2)$}
\llabel{ss-SU2}

The Lie group $SU(2)$ is the group of unitary unimodular $2\times 2$ complex matrices 
$$SU(2)=\Pg{\Pt{\ba{cc}
\al & \beta\\
-\overline{\beta} & \overline{\al}
\ea}\in\Mat(2,\C)\ |\ |\al|^2+|\beta|^2=1}.$$
It is compact and simply connected. The Lie algebra of $SU(2)$ is the algebra of antihermitian traceless $2\times 2$ complex matrices 
$$su(2)=\Pg{\Pt{\begin{array}{cc}
i\alpha & \beta\\
-\overline{\beta} & -i\alpha
\end{array}}\in\Mat(2,\C)\ |\ \alpha\in\R,\beta\in\C}.$$
A basis of $su(2)$ is $\Pg{p_1,p_2,k}$ where 
\bqn
p_1=\frac{1}{2}\Pt{\begin{array}{cc}
0 & 1\\
-1 & 0
\end{array}}\quad
p_2=\frac{1}{2}\Pt{\begin{array}{cc}
0 & i\\
i & 0
\end{array}}\quad
k=\frac{1}{2}\Pt{\begin{array}{cc}
i & 0\\
0 & -i
\end{array}},
\eqnn
whose commutation relations are $[p_1,p_2]=k\quad[p_2,k]=p_1\quad[k,p_1]=p_2$. 
Recall that for $su(2)$ we have $\Kil(X,Y)=4\Tr(XY)$ and, in particular, $\Kil(p_i,p_j)=-2\de_{ij}$.
The choice of the subspaces $$\k=\span{k}\qquad \p=\span{p_1,p_2}$$
provides a {\it Cartan decomposition} for $su(2)$.
Moreover, $\Pg{p_1,p_2}$ is an orthonormal frame for the inner product $<\cdot,\cdot>=-\frac{1}{2}\Kil(\cdot,\cdot)$ restricted to $\p$.

Defining $\Delta(g)=g\p$ and ${\mathbf g}_g(v_1,v_2)=< g^{-1}v_1,g^{-1}v_2>$, we have that $(SU(2),\Delta,{\mathbf g})$ is a $\k\oplus\p$ sub-Riemannian manifold (for details see Appendix \ref{s-k+p}). The sub-Riemannian manifold $SU(2)$ is thus endowed with the standard definition of sub-Riemannian length and distance.

\subsubsection{Expression of geodesics}

$\k\oplus\p$ manifolds have very special properties: there are no strict abnormal minimizers, hence all the geodesics starting from a point $g_0$ are parametrized by the initial covector. The explicit expression of a geodesic starting from $g_0$ is
\bqn g(t)=g_0e^{(A_k+A_p)t}e^{-A_k t}.
\eqnl{e-k+p}
with $A_k\in\k,\ A_p\in\p$. The geodesic is parameterized by arclength when $<A_p,A_p>=1$. This condition defines a cylinder $\Lam_{g_0}\subset T^*_{g_0}M$.

For the $SU(2)$ case, we compute the geodesics starting from the identity using  formula \r{e-k+p}. Consider an initial covector $\lam=\lam(\th,c)=\cos(\th)p_1+\sin(\th)p_2+ck\in \Lambda_\Id$. The corresponding exponential map is
\bqn
\EXP(\th,c,t)&:=&\EXP(\lam(\th,c),t)=e^{(\cos(\th)p_1+\sin(\th)p_2+ck)t} e^{-ck t}=\Pt{\ba{cc}\al&\beta\\-\overline{\beta}&\overline{\al}\ea}
\eqnn
with
\bqn
\Pt{\ba{c} \al\\\beta\ea}=\Pt{\begin{array}{c}
\frac{c\sin(\frac{ct}{2})\sin (\sqrt{1 + c^2}\frac{t}{2})}{\sqrt{1 + c^2}} + \cos(\frac{ct}{2})\cos(\sqrt{1 + c^2}\frac{t}{2})+i\Pt{\frac{c\cos(\frac{ct}{2})\sin (\sqrt{1 + c^2}\frac{t}{2})}{\sqrt{1 + c^2}} - \sin(\frac{ct}{2})\cos(\sqrt{1 + c^2}\frac{t}{2})}\\
\frac{\sin(\sqrt{1 + c^2}\frac{t}{2})}{\sqrt{1 + c^2}}
\Pt{\cos(\frac{ct}{2}+\th) + i\sin(\frac{ct}{2}+\th)}
\end{array}}.\eqnn

Another property of $\k\oplus\p$ manifolds is the existence of a solution of the minimal length problem (see Remark \ref{rem-minesiste}), i.e. for each pair $g,h\in M$ there exists a trajectory $\ga$ steering $g$ to $h$ and minimizing the length: $l(\ga)=d(g,h)$. Recall that this minimizing trajectory is a geodesic.

\subsubsection{The cut locus and distance for $SU(2)$}

In this section we recall the formula of the cut locus for the $\k\oplus\p$ manifold $SU(2)$ and of the sub-Riemannian distance. Both the results are proved in \cite{nostro-gruppi}.
\bt
The cut locus for the $\k\oplus \p$ problem on $SU(2)$ is $$K_\Id=e^\k\setminus\Id=\Pg{e^{ck}\ |\ c\in(0,4\pi)}.$$
\et

The cut locus $K_\Id$ is topologically equivalent to a circle $S^1$ without a point, the starting point $\Id$.

\bt
\llabel{t-distance}
Let $g=\matSU\in SU(2)$. Consider the sub-Riemannian distance from $\Id$ defined by the $\k\oplus\p$ structure on $SU(2)$. It holds
\bqn
d\Pt{\Id,g}=\begin{graffa3}
2\sqrt{\arg(\al)\Pt{2\pi-\arg(\al)}}& \mbox{~if~} & \beta=0\\
\psi(\al) &\mbox{~if~} & \beta\neq 0,
\end{graffa3}
\eqnn
where $\arg\Pt{\al}\in[0,2\pi]$ and $\psi(\al)=t$ where $t$ is the unique solution of the system
\bqn
\begin{cases}
-\frac{ct}{2}+\arctan\Pt{\frac{c}{\sqrt{1+c^2}}\tan\Pt{\frac{\sqrt{1+c^2}t}{2}}}=\arg(\al)\\
\frac{\sin\Pt{\frac{\sqrt{1+c^2}t}{2}}}{\sqrt{1+c^2}}=\sqrt{1-|\al|^2}\\
t\in\Pt{0,\frac{2\pi}{\sqrt{1+c^2}}}.
\end{cases}
\eqnn
\et

\brem Notice that $\forall~\al,\beta_1,\beta_2\in\C, |\beta_1|=|\beta_2|$ we have $d\Pt{\Id,\Pt{\ba{c}\al\\ \beta_1\ea}}=d\Pt{\Id,\Pt{\ba{c}\al\\ \beta_2\ea}}=d\Pt{\Id,\Pt{\ba{c}\overline{\al}\\\beta_1\ea}}$.
\erem

\subsection{A sub-Riemannian problem on $\Lqu$}

We define the 3-D sub-Riemannian manifold $\Lqu$ as a quotient of the $\k\oplus\p$ manifold $SU(2)$ defined above. $\Lqu$ is a lens space: for more details about lens spaces $L(p,q)$ see \cite{rolfsen}.
We prove that the quotient and the sub-Riemannian structure are compatible. We then compute the distance $\hat{d}(.,.)$ on $\Lqu$. 

\subsubsection{Definition of \Lqu}

We define coordinates $(a,b,c)$ on $SU(2)$ as follows:
\bp
Each $g\in SU(2)$ can be written as
\bqn
g=e^{-bp_2}e^{ap_1}e^{cp_2}
\eqnl{e-SU2coord}
with $a\in[0,\pi]$, $b\in[0,2\pi)$, $c\in\R/({4\pi})$. The value of $a$ is uniquely determined by $g$.

If $a\neq 0,~a\neq \pi$, then the values of $b$ and $c$ are uniquely determined by $g$.

If $a=0$, then the value of $c-b \mod 4\pi$ is uniquely determined by $g$.

If $a=\pi$, then the value of $c+b \mod 4\pi$ is uniquely determined by $g$.
\ep
\proof
Observe that $e^{-bp_2}e^{ap_1}e^{cp_2}=\matSU$ with
\bqn
\bcas
\al=\cos\Pt{\frac{a}{2}}\cos\Pt{\frac{c-b}{2}}+i\sin\Pt{\frac{a}{2}}\sin\Pt{\frac{c+b}{2}}\\
\beta=\sin\Pt{\frac{a}{2}}\cos\Pt{\frac{c+b}{2}}+i\cos\Pt{\frac{a}{2}}\sin\Pt{\frac{c-b}{2}}.
\ecas
\eqnn 
The result follows from direct computation.
\qed

In the following we will write $g=(a,b,c)$ for $g=e^{-bp_2}e^{ap_1}e^{cp_2}$, even though the vector $(a,b,c)$ is not uniquely determined by $g$.
\bp
\label{p-L41def}
Define the following equivalence relation on $SU(2)$: $g_1\approx g_2$ if there exist $(a_1,b_1,c_1),(a_2,b_2,c_2)$ such that $g_1=(a_1,b_1,c_1),~g_2=(a_2,b_2,c_2)$ and $\begin{cases}
a_1=a_2\\
b_1=b_2\\
c_1=c_2\mod \pi.
\end{cases}$

The 3-D manifold $SU(2)/\approx$ is the lens space $\Lqu$.
\ep
Before proving the proposition, we recall the standard definition of $L(p,q)$ with $p,q\in\Z\backslash\Pg{0}$ coprime. Let $$S^3=\Pg{(x_1,x_2)\in\C^2\ |\ |x_1|^2+|x_2|^2=1}$$ and define on $S^3$ the following equivalence relation: $(x_1,x_2)\sim(y_1,y_2)$ if there exists $\om\in\C$ $p$-th root of unity (i.e. $\om^p=1$) such that
$$\Pt{\ba{c}x_1\\x_2\ea}=\Pt{\ba{cc}
\om& 0\\
0&\om^q\ea}\Pt{\ba{c}y_1\\y_2\ea}.$$
The quotient manifold $S^3/\sim$ is the lens space $L(p,q)$.

{\bf Proof of Proposition \ref{p-L41def}.} Consider the isomorphism \mfunz{\phi}{SU(2)}{S^3}{\matSU}
{\Pt{\frac{\al+\beta}{\sqrt{2}},\frac{\Re{\al}+i\Im{\beta}-\Re{\beta}-i\Im{\al}}{\sqrt{2}}}.} 
A straightforward computation gives that $g_1,g_2\in SU(2)$ are equivalent ($g_1\approx g_2$) if and only if $\phi(g_1),\phi(g_2)\in S^3$ are equivalent $(\phi(g_1)\sim\phi(g_2))$ with respect to the equivalence relation defined by $p=4$, $q=1$. Hence the manifolds $SU(2)/_\approx$ and $S^3/_\sim$ are isomorphic. Thus $SU(2)/_\approx$ is a 3-D manifold, the lens space $\Lqu$.
\qed

\subsubsection{Sub-Riemannian structure on $\Lqu$}
We now define a sub-Riemannian structure on $\Lqu$ induced by the quotient map.
\bp
\llabel{p-projSU2Lpq}
The sub-Riemannian structure on $SU(2)$ given in \ref{ss-SU2} induces a 2-dim sub-Riemannian structure on $\Lqu$ via the quotient map~~\mmfunz{\Pi}{SU(2)}{\Lqu}{g}{[g],}
i.e.
\bi
\i the map $$\widetilde{\Delta}:[g]\mapsto\Pi_*\Pt{\Delta(h)}\subset T_{[g]}\Lqu\mbox{~with~}h\in[g]$$ is a 2-dim smooth distribution on $\Lqu$ that is Lie bracket generating;
\i $\tilde{\g}_{[g]}(v_*,w_*)=\Pa{v_*,w_*}_{[g]}:=\Pa{v,w}_h\mbox{~with~} h\in[g],~v,w\in T_hSU(2),~\Pi_*(v)=v_*,~\Pi_*(w)=w_*$ is a smooth positive definite scalar product on $\widetilde{\Delta}$.
\ei
\ep
\proof
The role of the map $\Pi$ and $\Pi_{*_{|_g}}$ is illustrated in the following diagram
$$\xymatrix{T_gSU(2) \ar[d] \ar[r]^{\Pi_{*_{|_g}}} & T_{[g]}\Lqu \ar[d] \\
SU(2) \ar[r]^\Pi &  \Lqu.}$$

The map $\Pi$ is a local diffeomorphism, thus $\Pi_{*_{|_g}}:T_gSU(2)\to T_{[g]}\Lqu$ is a linear isomorphism, hence $\Pi_{*_{|_g}}\Pt{\Delta(g)}$ is a 2-dim subspace of $T_{[g]}\Lqu$.

The following statements:
\bi
\i the distribution $\widetilde{\Delta}([g])$ is well defined, i.e. $\forall~ h_1,h_2\in[g]$ we have $$\Pi_{*_{|_{h_1}}}\Pt{\Delta(h_1)}=\Pi_{*_{|_{h_2}}}\Pt{\Delta(h_2)};$$
\i the positive definite scalar product $\Pa{v_*,w_*}_{[g]}$ is well defined, i.e. $\forall~ h_1,h_2\in[g],~~v_1,w_1\in T_{h_1}SU(2),~~v_2,w_2\in T_{h_2}SU(2)$ such that $\Pi_{*_{|_{h_1}}}(v_1)=\Pi_{*_{|_{h_2}}}(v_2)$ and $\Pi_{*_{|_{h_1}}}(w_1)=\Pi_{*_{|_{h_2}}}(w_2)$ we have $\Pa{v_1,w_1}_{h_1}=\Pa{v_2,w_2}_{h_2};$\ei
are consequences of the following lemma:
\bl Let $h_1,h_2\in[g]$ with $h_1=(a,b,c)$ and $h_2=(a,b,c+k\pi)$, $k\in\Z$. The map \mmfunz{\phi}{\p}{\p}{m=m_1p_1+ m_2p_2}{(-1)^k m_1p_1+m_2p_2} is bijective and it is an isometry w.r.t. the positive definite scalar product $\Pa{~,~}$. The following equality holds $\forall~m\in\p$ $$\frac{d}{dt}_{|_{t=0}}\Pq{h_1e^{tm}}=\frac{d}{dt}_{|_{t=0}}\Pq{h_2e^{t\phi(m)}}.$$
\el
\proof
A direct computation gives that $\phi$ is bijective and that it is an isometry.

Consider now $h_1e^{tm}$ and $h_2 e^{tn}$ with $m=m_1p_1+m_2p_2\in\p$, $n=n_1p_1+n_2p_2\in\p$. We have
\bqn
h_1e^{tm}&=&e^{-bp_2}e^{ap_2}e^{cp_2}e^{tm}=
\Pt{\ba{cc}\al(a,b,c,t,m)&\beta(a,b,c,t,m)\\-\overline{\beta}(a,b,c,t,m)&\overline{\al}(a,b,c,t,m)\ea}_,\eqnn
with
\bqn
\al(a,b,c,t,m)&=&K_0(t,m)\al(a,b,c,0,0)+(-m_1+im_2)K_1(t,m)\beta(a,b,c,0,0),\nn
\beta(a,b,c,t,m)&=&K_0(t,m)\beta(a,b,c,0,0)+(m_1+im_2)K_1(t,m)\beta(a,b,c,0,0),\nn
K_0(t,m)&:=&\cos\Pt{\frac{t\sqrt{m_1^2+m_2^2}}{2}},~~~~~~
K_1(t,m):=\frac{\sin\Pt{\frac{t\sqrt{m_1^2+m_2^2}}{2}}}{\sqrt{m_1^2+m_2^2}}.
\eqnn

Observe that $h_2e^{tn}=e^{-bp_2}e^{ap_2}e^{cp_2+k\pi p_2}e^{tn}\simeq\Pt{\ba{cc}\al_2&\beta_2\\
\-\bar\beta_2&\bar\al_2\ea}$ with 
\bqn
\al_2&=&K_0(t,n)\al(a,b,c,0,0)+((-1)^{k+1}n_1+in_2)K_1(t,n)\beta(a,b,c,0,0),\nn
\beta_2&=&K_0(t,n)\beta(a,b,c,0,0)+((-1)^k n_1+in_2)K_1(t,n)\beta(a,b,c,0,0),
\eqnn
hence $\Pq{h_1e^{tm}}=\Pq{h_2e^{tn}}$ if both the following conditions are satisfied: $\bcas
K_0(t,m)=K_0(t,n)\\
K_1(t,m)=K_1(t,n)\\
-m_1+im_2=(-1)^{k+1}n_1+in_2\\
m_1+im_2=(-1)^{k}n_1+in_2.
\ecas$

These are verified when $n=\phi(m)$. Thus $\frac{d}{dt}_{|_{t=0}}\Pq{h_1e^{tm}}=\frac{d}{dt}_{|_{t=0}}\Pq{h_2e^{t\phi(m)}}$.
\eproof\\

Since $\Pi$ is a local diffeomorphism, $\forall~g\in SU(2)~~~\exists~B(g)$ such that the map $\Pi_{*_{|_{B(g)}}}:T_{B(g)}SU(2)\to T_{B([g])}\Lqu$ is a diffeomorphism, thus $\widetilde{\Delta}$ is smooth, Lie bracket generating, and $\Pa{v_*,w_*}_{[g]}$ is smooth as a function of $[g]$.
\eproof
\brem\llabel{r-Lqudiversi}
Observe that other definitions of a sub-Riemannian structure on $\Lqu$ induced by the one on $SU(2)$ are possible: for example, in \cite{nostro-gruppi} we defined the following identification rule on $SU(2)$:
\bqn
\matSU\sim \Pt{\ba{cc}\al '&\beta '\\-\bar\beta '& \bar \al '\ea}
\mbox{~if exists $\om\in\C$ such that $\om^p=1$ and~~~}\Pt{\ba{c}\al '\\\beta ' \ea}=
\Pt{\ba{cc}\om & 0\\ 0 & \om^q\ea}\Pt{\ba{c}\al \\\beta\ea}
\eqnl{e-Lqusbagliato}
and observed that the sub-Riemannian structure defined on $SU(2)$ in \ref{ss-SU2} induces a sub-Riemannian structure on the manifold $L(p,q):=SU(2)/\sim$.

The sub-Riemannian structures are locally isomorphic, but globally the structure may vary: indeed, in the following we will prove that the cut locus for the structure defined in \ref{p-projSU2Lpq} is globally different from the one defined by \r{e-Lqusbagliato}.
\erem
\brem Notice that in this case the push-forward of a left invariant vector field is not a vector field, because the projections $\Pi_{*_{|_g}}(gp),\ \Pi_{*_{|_h}}(hp)$ from different points $g,h$ such that $[g]=[h]$ do not coincide. Nevertheless, the projections of the whole distribution $\Pi_{*_{|_g}}(\Delta(g)),\Pi_{*_{|_h}}(\Delta(h))$ coincide; then the sub-Riemannian structure on $\Lqu$ is well defined even though the projections of the left invariant field is not well defined.
\erem

We have standard definitions of sub-Riemannian length and distance on $\Lqu$, as presented in Appendix \ref{a-SubRiem}, see \r{e-lunghezza} and \r{e-dipoi}. We indicate them respectively with $\hat{l}$ and $\hat{d}$.

Observe that the geodesics for the sub-Riemannian manifold $\Lqu$ are the projections of the geodesics for $SU(2)$, due to the fact that the projection $\Pi$ is a local isometry. In general we have the following
\bp
Let $\fz{\Gamma}{[0,T]}{SU(2)}$ be a smooth curve. Define its projection $\fz{\underline{\Ga}}{[0,T]}{\Lqu}$ by $\underline{\Ga}(t):=\Pi(\Ga(t))$. Then $\underline{\Ga}$ is a smooth curve and its length coincide with the length of $\Ga$:
$$\hat{l}(\underline{\Ga})=l(\Ga).$$
\ep

Recall now the definition of the {\bf lift} of a curve in our case and a standard result about its optimality:
\bdeff
Let $\fz{\ga}{[0,T]}{\Lqu}$ be a smooth curve in $\Lqu$ with $\ga(0)=[g]$. Fix $g_1\in SU(2)$ such that $g_1\in[g]$. The lift of $\ga$ starting from $g_1$ is the unique smooth curve $\fz{\bar\ga}{[0,T]}{SU(2)}$ satisfying $\bar\ga(0)=g_1$ and $\Pi\Pt{\bar{\ga}(t)}=\ga(t)$ for any $t\in[0,T]$.

The length of the lift $\bar\ga$ in $SU(2)$ coincide with the length of the curve $\ga$ in $\Lqu$:
$$l(\bar\ga)=\hat{l}(\ga).$$
\edeff
\bp
Let $\fz{\ga}{[0,T]}{\Lqu}$ be an optimal trajectory steering $[g]$ to $[h]$. Fix $g_1\in SU(2)$ with $g_1\in[g]$. Then the lift $\fz{\bar\ga}{[0,T]}{SU(2)}$ starting from $g_1$ is an optimal trajectory steering $g_1$ to $\bar\ga(T)\in[h]$.
\ep

We end this section giving an expression to compute the sub-Riemannian distance on $\Lqu$.
\bp
\label{t-distLqu}
The sub-Riemannian distance $\hat{d}$ on $\Lqu$ can be computed via the sub-Riemannian distance $d$ on $SU(2)$ as follows:
\bqn
\hat{d}\Pt{[g],[h]}=min_{h_1\in[h]}\Pg{d(g_1,h_1)}
\eqnn
independently on the choice of $g_1\in[g]$.
\ep
\proof
Consider a trajectory $\ga$ in $\Lqu$ that is solution of the minimal length problem between $[g]$ and $[h]$. Define its lift $\bar\ga(.)$ starting from a fixed $g_1\in[g]$. It steers $g_1$ to a given $h_1\in[h]$ and it is optimal, thus $\hat{d}([g],[h])=d(g_1,h_1)$.

We prove now that for any $h_2\in[h]$ we have $d(g_1,h_1)\leq d(g_1,h_2)$. By contradiction, assume that $\exists\ h_2 \in [h]$ such that $d(g_1,h_2)<d(g_1,h_1)$, thus there exists a smooth curve $\fz{\Ga}{[0,T]}{SU(2)}$ steering $g_1$ to $h_2$ with $l(\Ga)<d(g_1,h_1)$. Consider its projection $\underline{\Ga}$: it steers $[g]$ to $[h]$ and it satisfies $\hat{l}(\underline{\Ga})<\hat{d}([g],[h])$. Contradiction.
\qed

\subsection{Isometry between $\Lqu$ and $PTS^2$}
This section is fundamental in the treatment of the mechanical problem presented in Section \ref{s-mecc}: we state that the problem \r{OCsfera} on $PTS^2$ is the problem of minimal length on the sub-Riemannian manifold $\Lqu$.

We start writing the expression in coordinates of the projection $S^3\rightarrow S^2$,  the well-known Hopf map.
\bp
Consider the Hopf map \mmfunz{\psi}{SU(2)\simeq S^3}{S^2}{\matSU}{
\Pt{\ba{c}
2\Pt{\Re{\al}\Re{\beta}+\Im{\al}\Im{\beta}}\\
2\Pt{\Re{\al}\Im{\al}-\Re{\beta}\Im{\beta}}\\
\Re{\al}^2+\Im{\beta}^2-\Re{\beta}^2 -\Im{\al}^2.
\ea
}}
Consider the coordinates $(a,b,c)$ on $SU(2)$ defined in \r{e-SU2coord} and the spherical coordinates $(a,b)$ on $S^2$ defined in \r{sfer}: then the Hopf map expression in coordinates is
\mmfunz{\psi}{SU(2)}{S^2}{(a,b,c)}{(a,b).}
\ep

\brem
Observe that both the $(a,b,c)$ coordinates on $SU(2)$ and the spherical coordinates $(a,b)$ are not well defined in $a=0$ and $a=\pi$. Nevertheless, the map $\phi$ is well defined and has this expression in coordinates also in the degenerate cases.
\erem

The following theorem is the central result of this section: here we prove that our mechanical problem \r{OCsfera} coincide with the sub-Riemannian problem on $\Lqu$.

\bt
The lens space $\Lqu$ is diffeomorphic to  $PTS^2$, the bundle of directions of the sphere. The coordinates $(a,b,c)$ on $\Lqu$ coincide with the coordinates $(a,b,-c)$ on $PTS^2$. The diffeomorphism maps horizontal curves on $\Lqu$ to horizontal curves on $PTS^2$ and vice versa.

The sub-Riemannian length of an horizontal curve $(a(.),b(.),c(.))$ in $\Lqu$ is equal to the cost $\J[\Ga]$ of the curve $\Ga(.)=(a(.),b(.),-c(.))$ in $PTS^2$.
\et
\proof
Consider any horizontal smooth curve $\fz{\Ga}{(-\eps,\eps)}{SU(2)}$ with $\Ga(0)=(a,b,c)$. We have $\Ga(t)=e^{-bp_2}e^{ap_1}e^{cp_2}e^{t\eta + o(t)}$ with $\eta\in\p$ due to the fact that $\Ga$ is horizontal. Assume $\eta=\eta_1p_1+\eta_2p_2$  satisfying $\eta_1\neq 0$. The reason of this condition will be clear in the following.

Consider its projection $\fz{\ga=\psi\circ\Ga}{(-\eps,\eps)}{S^2}$: we have by computation
\bqn\ga(t)=\Pt{\ba{c}
\sin(a)\cos(b)\\
\sin(a)\sin(b)\\
\cos(a)
\ea
}+t\eta_1\Pt{\ba{c}
\cos(a)\cos(b)\cos(c)+\sin(b)\sin(c)\\
\cos(a)\sin(b)\cos(c)-\cos(b)\sin(c)\\
-\sin(a)\cos(c)\ea}+o(t).
\eqnn

Due to condition $\eta_1\neq 0$, $\ga$ is smooth and differentiable in $0$; the direction of its tangent vector in $0$ depends only on $\Ga(0)=(a,b,c)$. Its lift is $(a(.),b(.),\xi(.))$ with $\xi(t)=\arccos\Pt{\frac{\dot{b}\sin(a)}{\dot{a}}}$: an explicit computation gives $\xi(t)=-\arctan(\tan(c))$. Thus the function \mfunz{f}{SU(2)}{PTS^2}{(a,b,c)}{(a,b,-\arctan(\tan(c)))} maps horizontal smooth curves of $SU(2)$ with $\eta_1\neq 0$ to horizontal smooth curves of $PTS^2$.

The case $\eta_1=0$ can be studied as a limit case: in this case $\dot{a}=\dot{b}=0$, thus $f$ maps horizontal curves of $SU(2)$ to the curves of $PTS^2$ representing rotations on a point (see Remark \ref{rem-PTS2fixedpoint}). In this sense $f$ maps all horizontal smooth curves of $SU(2)$ to horizontal smooth curves of $PTS^2$.

Consider now $\Lqu=SU(2)/_\approx$: the map $f$ pass to the quotient and we have the explicit expression \mmfunz{f_\approx}{\Lqu}{PTS^2}{[a,b,c]}{(a,b,-c).} Observe that $f_\approx$ is a diffeomorphism. Moreover, it maps horizontal smooth curves of $\Lqu$ to horizontal smooth curves of $PTS^2$. Then the differential $df_\approx$ is a linear isomorphism $\fz{df_\approx}{\widetilde{\Delta}([g])\subset\Lqu}{\Delta(a,b,-c)\subset PTS^2}$ that is explicitly $df_\approx(\eta_1 \Pi_{*_{|_g}}(gp_1)+\eta_2 \Pi_{*_{|_g}}(gp_2))=\eta_1F_1(a,b,-c) -\eta_2F_2(a,b,-c)$ with $g=(a,b,c)\in SU(2)$. Then we have that the length of a curve $\fz{\ga}{[0,T]}{\Lqu}$ satisfies $$\hat{l}(a(.),b(.),c(.))=\int_0^T\sqrt{\eta_1^2+\eta_2^2}=J[a(.),b(.)].$$
\qed

\section{Proof of the main result}
\llabel{s-solution}

We give now the explicit computation of the cut locus for the problem {\bf (P)} on $S^2$. For this specific problem the first conjugate locus (see Appendix \ref{a-conj}) is completely contained in the cut locus since this is the case for the lifted problem on $SU(2)$ (see \cite[Sections 3.1.3, 5.1]{nostro-gruppi}). Hence the cut locus is the set of points where geodesics lose optimality (see Remark \ref{rem-agra}).

\bt
\llabel{t-L41cut}
The cut locus $K_{[\Id]}$ for the sub-Riemannian problem on $\Lqu$ starting from the identity is a stratification $K_{[\Id]}=K^{loc}\cup K^{sym}$ with
\bqn
K^{loc}&=&[K_\Id]\backslash[\Id]=[e^\k]\backslash[\Id],\nn
 K^{sym}&=&\left\{\Pq{g}\in\Lqu\ |\mbox{~for~a~fixed~}g^*\in[g]~\exists~x_1,x_2\in [\Id], x_1\neq x_2 \mbox{ such that }\right.\nn
&&\hspace{5cm}\left.\forall x^*\in[\Id]~~d(g^*,x_1)=d(g^*,x_2)\leq d(g^*,x^*)\right\}.
\eqnn
where $d(.,.)$ is the sub-Riemannian distance on $SU(2)$ given in Theorem \ref{t-distance}.
\et
\proof
We first prove that $K^{loc}$ lies in the cut locus. Consider $[g]\in K^{loc}$: due to the definition of $[g]$, there exists $g_1\in[g]$ such that $g_1=e^{nk}=\Pt{\ba{cc}
\al & 0\\
0& \bar{\al}\ea}$ with $|\al|=1$. Hence $[g]=\Pg{g_1,-g_1,g_2,-g_2}$ with $g_2=\Pt{\ba{cc} 0 & i\al\\
i\al & 0\ea}$. Recall that $\hat{d}([\Id],[g])=\min_{g^*\in[g]}d(\Id,g^*)$. Applying Theorem \ref{t-distance}, we have $d(\Id,g_1)=2\sqrt{\arg(\al)(2\pi-\arg(\al))}$, $d(\Id,-g_1)=\bcas
2\sqrt{(\arg(\al)+\pi)(\pi-\arg(\al))}\mbox{~~~if $\al\in[0,\pi]$}\\
2\sqrt{(\arg(\al)-\pi)(3\pi-\arg(\al))}\mbox{~~~if $\al\in[\pi,2\pi)$}
\ecas$, $d(\Id,g_2)=d(\Id,-g_2)=\pi$. A straightforward computation gives that the minimum is attained for $d(\Id,g_1)$ or $d(\Id,-g_1)$. We assume without loss of generality that we have $\hat{d}([\Id],[g])=d(\Id,g_1)$. Recall that $g_1\in K_\Id$, thus there exist $\ga_1,\ga_2$ two different optimal trajectories on $SU(2)$ steering $\Id$ to $g_1$. Thus their projections $\underline{\ga_1},\underline{\ga_2}$ are trajectories on $\Lqu$ steering $[\Id]$ to $[g]$ and satisfying $\hat{l}(\underline{\ga_i})=l(\ga_i)=d(\Id,g_1)=\hat{d}([\Id],[g])$ for $i=1,2$; hence they are both optimal. Observe that $\ga_1,\ga_2$ are distinct in a \nei\ of $\Id$, thus $\underline{\ga_1},\underline{\ga_2}$ are distinct in a \nei\ of $[\Id]$. Thus $[g]$ is a cut point, due to the fact that it is reached by two different optimal trajectories $\underline{\ga_1},\underline{\ga_2}$.

We prove now that $K^{sym}$ lies in the cut locus. Consider a point $[g]\in K^{sym}$ and a fixed $g^*\in[g]$. Recall that $d(g^*,x_1)=d(g^*,x_2)$. Due to the existence of a solution to the minimal length problem, there exist $\Ga_1$ and $\Ga_2$ geodesics in $SU(2)$ with the following property: $\Ga_i$ steers $x_i\in[\Id]$ to $g^*\in[g]$ and it is optimal. Observe that $l(\Ga_1)=l(\Ga_2)$ and that $\Ga_1$ and $\Ga_2$ don't coincide in a \nei\ of $g^*$ due to the fact that they are two different optimal geodesics. Now consider the two projections $\underline{\Ga_1}$ and $\underline{\Ga_2}$: they are trajectories in $\Lqu$ steering $[\Id]$ to $[g]$ and their length is $\hat{l}(\underline{\Ga_i})=d(g^*,x_i)=\min_{x_j\in[\Id]}d(g^*,x_j)$ due to the fact that 
$d(g^*,x_1)=d(g^*,x_2)\leq d(g^*,x^*).$ Hence $\underline{\Ga_1}$ and $\underline{\Ga_2}$ are optimal geodesics. Notice that they are distinct in a \nei\ of $[g]$ because they are projections of two trajectories that are distinct in a \nei\ of $g^*$. Observe that the choice of $g^*\in[g]$ is not relevant, as a consequence of Theorem \ref{t-distLqu}.

We end proving that there is no cut point outside $K^{loc}\cup K^{sym}$. By contradiction, let $[g]\in \Lqu\backslash(K^{loc}\cup K^{sym}\cup\Pg{[\Id]})$ be a cut point: thus there are $\ga_1$ and $\ga_2$ two different optimal geodesics of $\Lqu$ steering $[Id]$ to $[g]$. Consider their two lifts $\overline{\ga_1}$ and $\overline{\ga_2}$ starting from $\Id\in SU(2)$. They reach respectively $g_1,g_2\in[g]$ and they are optimal in $SU(2)$. We consider two distinct cases:
\bi
\i $g_1=g_2$: we have two distinct optimal trajectories $\overline{\ga_1}$ and $\overline{\ga_2}$ steering $\Id$ to the point $g_1$. Hence $g_1\in K_\Id$. Observe that $g_1\not\in[\Id]$ otherwise $[g]=[\Id]$. Thus $g_1\in K_\Id\backslash[\Id]$, hence $[g]\in K^{loc}$. Contradiction.

\i $g_1\neq g_2$: consider the two lifts of $\ga_1$ and $\ga_2$ reaching both the point $g_1$ and call them $\widetilde{\ga_1}$ and $\widetilde{\ga_2}$ respectively. Observe that $\widetilde{\ga_1}=\overline{\ga_1}$, hence $\widetilde{\ga_1}(0)=\Id$. Instead, $\widetilde{\ga_2}(0)=x_2\neq \Id$: we prove it by contradiction. If $\widetilde{\ga_2}(0)=\Id$, then $\widetilde{\ga_2}=\overline{\ga_2}$ by the uniqueness of the lift. Hence $\widetilde{\ga_2}(T)=g_2\neq g_1$, that contradicts the definition of $\widetilde{\ga_2}$. Contradiction. Observe that $x_2\in [\Id]$ due to the fact that $\ga_2(0)=[\Id]$.

We have that both $\widetilde{\ga_1}$ and $\widetilde{\ga_2}$ are optimal, hence $d(\Id,g_1)=l(\widetilde{\ga_1})=\hat{l}(\ga_1)=
\hat{l}(\ga_2)=l(\widetilde{\ga_2})=d(x_2,g_1)$ with $x_2\in [\Id]$ and $x_2\neq\Id$. Observe that, choosing any $x^*\in[\Id]$, we have $d(\Id,g_1)\leq d(x^*,g_1)$ due to the optimality of $\ga_1$. Thus $[g]\in K^{sym}$. Contradiction.
\ei
\eproof

We now show some pictures of the projections on $S^2$ of optimal geodesics on $\Lqu$, i.e. the optimal trajectories for our problem.

In Figure \ref{fig:cutlocale} some optimal trajectories with cusps are shown: they meet in the same point with the same direction, where they lose optimality. This point with direction lies in the local cut locus $K^{loc}$.

In Figure \ref{fig:cutlisci} we show smooth optimal trajectories, that are indeed half-circles: the points where they lose optimality lie in the intersection $K^{loc}\cap K^{a}$.

In Figure \ref{fig:cutentrambi} we show optimal trajectories (both smooth and with cusps) meeting in one of the two points of intersections of the three strata $K^{loc}\cap K^{a}\cap K^{b}$.

\immagine[7]{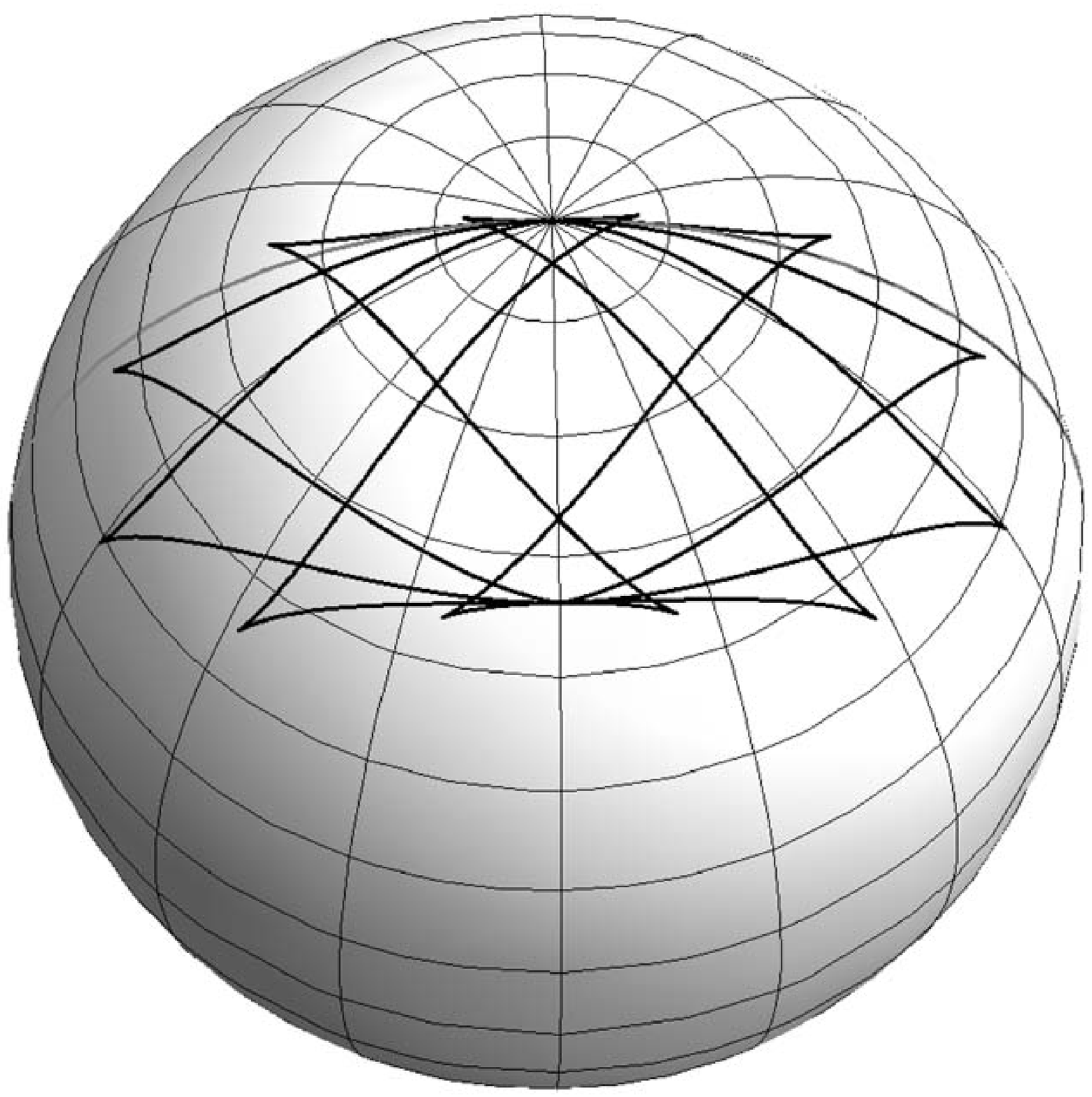}{Projections of optimal trajectories with cusps.}
\immagine[7]{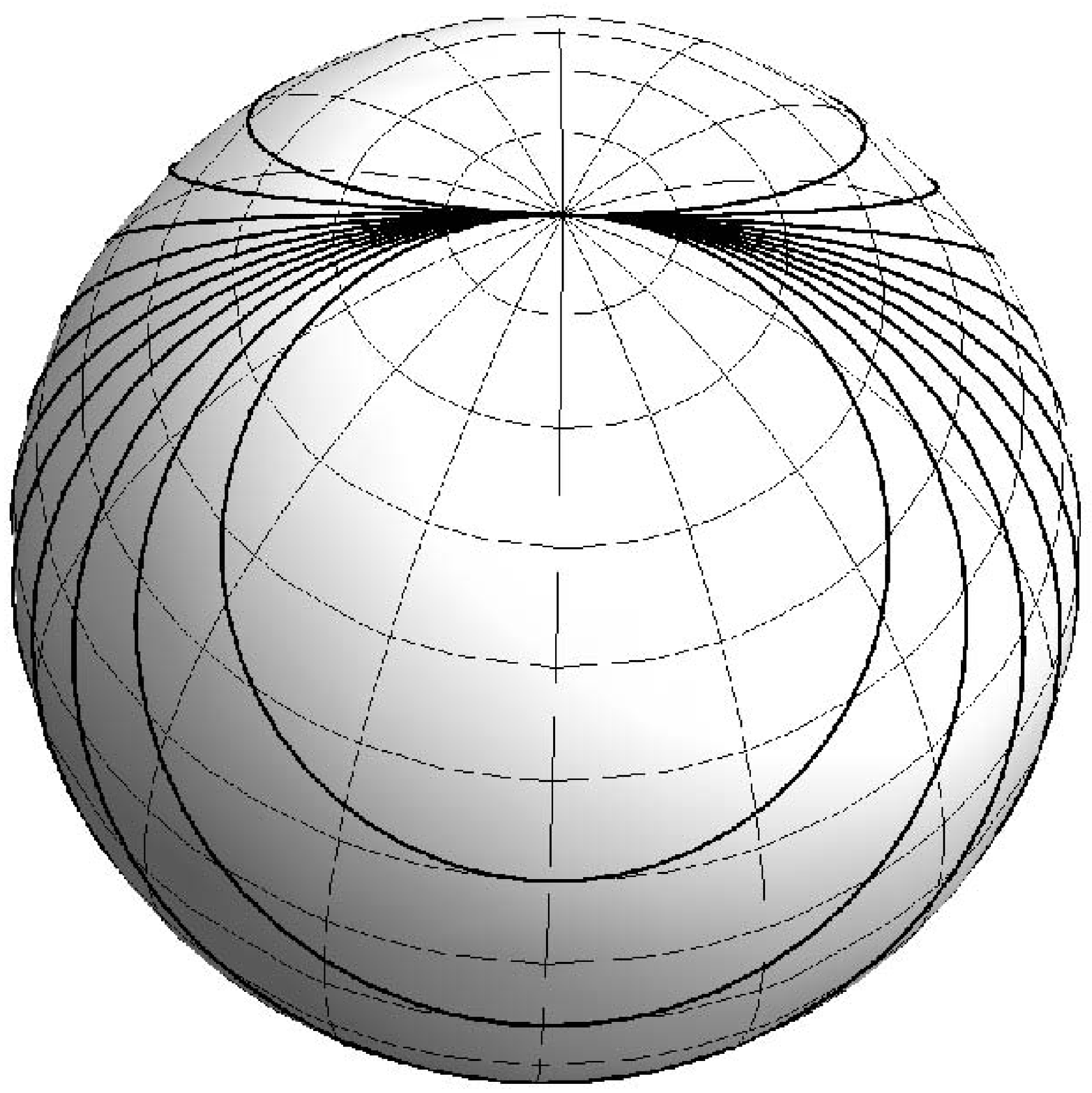}{Projections of smooth optimal trajectories.}
\immagine[7]{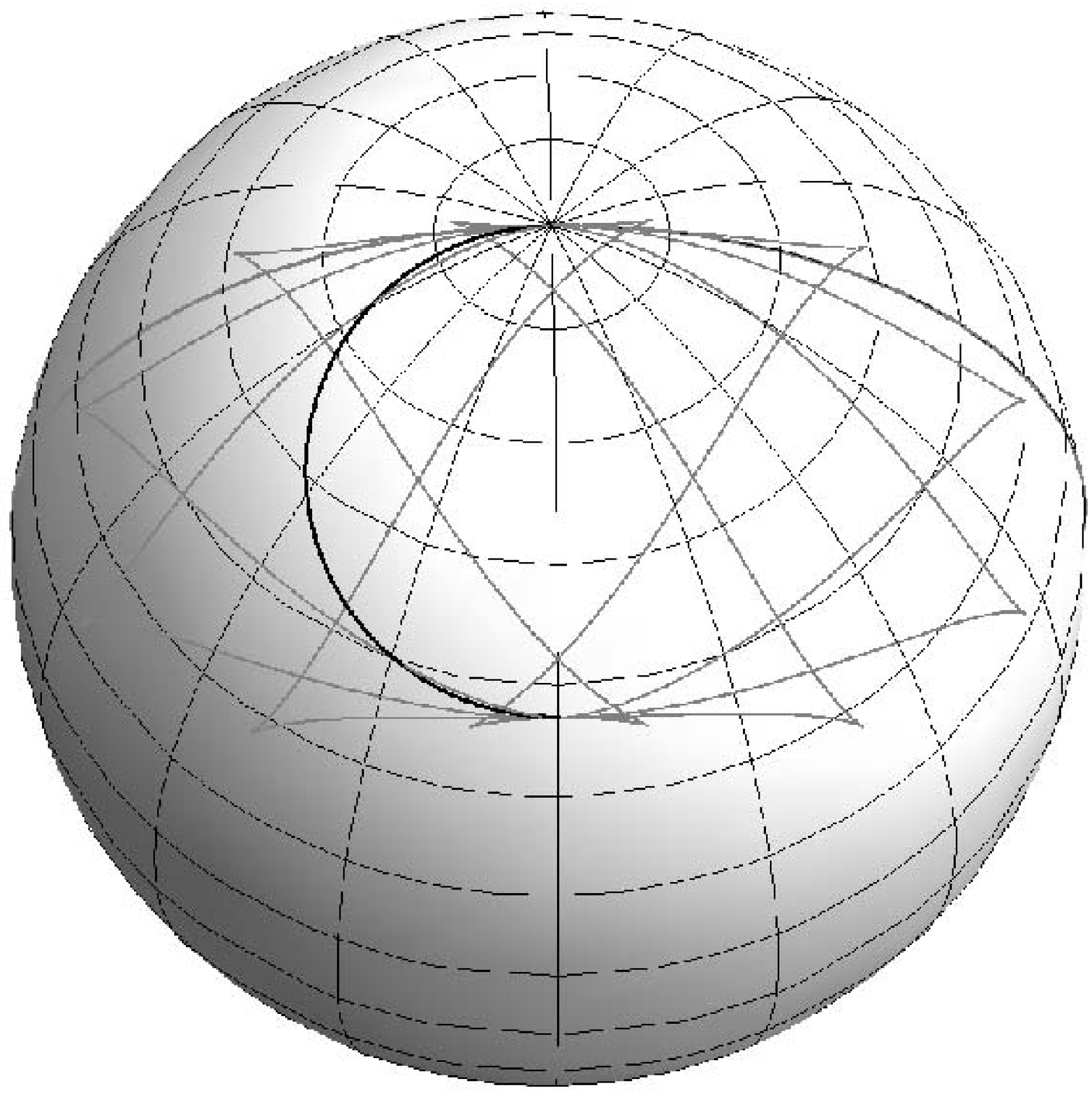}{Projections of optimal trajectories (both smooth and with cusps) meeting in $K^{loc}\cap K^{a}\cap K^{b}$.}

\appendix
\section{Basic results in sub-Riemannian geometry}
\llabel{a-SubRiem}

In this Appendix we recall some definitions and results about sub-Riemannian geometry, in particular about $\k\oplus\p$ manifolds. For a deeper presentation of sub-Riemannian geometry see e.g. \cite{bellaiche,gromov}.

\subsection{Sub-Riemannian manifold}
A $(n,m)$-sub-Riemannian manifold is a triple $(M,\Delta,{\mathbf g})$, 
where
\bi
\i $M$ is a connected smooth manifold of dimension $n$;
\i $\Delta$ is a Lie bracket generating smooth distribution of constant rank $m<n$, i.e. $\Delta$ is a smooth map that associates to $q\in M$  a $m$-dim subspace $\Delta(q)$ of $T_qM$, and $\forall~q\in M$ we have
\bqn\llabel{Hor}\span{[f_1,[\ldots[f_{k-1},f_k]\ldots]](q)~|~f_i\in\mathrm{Vec}(M)\mbox{~and~}f_i(p)\in\Delta(p)~\forall~p\in M}=T_qM.\eqn
Here $Vec(M)$ denotes the set of smooth vector fields on $M$.
\i ${\mathbf g}_q$ is a Riemannian metric on $\Delta(q)$, that is smooth 
as function of $q$.
\ei

The Lie bracket generating condition \r{Hor} is also known as H\"ormander condition.

A Lipschitz continuous curve $\ga:[0,T]\to M$ is said to be \b{horizontal} if 
$\dot\ga(t)\in\Delta(\ga(t))$ for almost every $t\in[0,T]$.
Given an horizontal curve $\ga:[0,T]\to M$, the {\it length of $\ga$} is
\bqn
l(\ga)=\int_0^T \sqrt{ \g_{\ga(t)} (\dot \ga(t),\dot \ga(t))}~dt.
\eqnl{e-lunghezza}
The {\it distance} induced by the sub-Riemannian structure on $M$ is the 
function
\bqn
d(q_0,q_1)=\inf \{l(\ga)\mid \ga(0)=q_0,\ga(T)=q_1, \ga\ \mathrm{horizontal}\}.
\eqnl{e-dipoi}

The \hp\ of connectedness of M and the Lie bracket generating assumption for the distribution guarantee the finiteness and the continuity of $d(\cdot,\cdot)$ with respect to the topology of $M$ (Chow's Theorem, see for instance \cite{agra-book}). The function $d(\cdot,\cdot)$ is called the Carnot-Charateodory distance and gives to $M$ the structure of metric space (see \cite{bellaiche,gromov}).

It is a standard fact that $l(\ga)$ is invariant under reparameterization of the curve $\ga$.
Moreover, if an admissible curve $\ga$ minimizes the so-called {\it energy functional}
$$ E(\ga)=\int_0^T {\g}_{\ga(t)}(\dot \ga(t),\dot \ga(t))~dt. $$
with $T$ fixed (and fixed initial and final point), then $v=\sqrt{\g_{\ga(t)}(\dot \ga(t),\dot \ga(t))}$
is constant and $\ga$ is also a minimizer of $l(\cdot)$.
On the other side a minimizer $\ga$ of $l(\cdot)$ such that  $v$ is constant is a minimizer of $E(\cdot)$ with $T=l(\ga)/v$.

A {\it geodesic} for  the sub-Riemannian manifold  is a curve $\ga:[0,T]\to M$ such that for every sufficiently small interval $[t_1,t_2]\subset [0,T]$, $\ga_{|_{[t_1,t_2]}}$ is a minimizer of $E(\cdot)$.
A geodesic for which $\g_{\ga(t)}(\dot \ga(t),\dot \ga(t))$  is (constantly) equal to one is said to be parameterized by arclength.

Locally, the pair $(\Delta,{\mathbf g})$ can be given by assigning a set of $m$ smooth vector fields that are orthonormal for ${\mathbf g}$, i.e.  
\bqn
\Delta(q)=\span{F_1(q),\dots,F_m(q)}, ~~~{\mathbf g}_q(F_i(q),F_j(q))=\delta_{ij}.
\eqnl{trivializable}
When  $(\Delta,{\mathbf g})$ can be defined as in \r{trivializable} by $m$ vector fields defined globally, we say that the sub-Riemannian manifold is {\it trivializable}. 


Given a $(n,m)$- trivializable sub-Riemannian manifold, the problem of finding a curve minimizing the energy between two fixed points  $q_0,q_1\in M$ is
naturally formulated as the optimal control problem 
\bqn
\dot q=\sum_{i=1}^m u_i F_i(q)\,,~~~u_i\in\R\,,
~~~\int_0^T
\sum_{i=1}^m u_i^2(t)~dt\to\min,~~q(0)=q_0,~~~q(T)=q_1.
\eqnl{sopra}

It is a standard fact that this optimal control problem is equivalent to the minimum time problem with controls $u_1,\ldots, u_m$ satisfying $u_1^2+\cdots+u_m^2\leq 1$.

When the manifold is analytic and the orthonormal frame can be assigned through $m$ analytic vector fields, we say that the sub-Riemannian manifold is {\it analytic}.

The $\k\oplus\p$ manifold presented below are trivializable and analytic since they are given in terms of left-invariant vector fields on Lie groups.
\subsection{First order necessary conditions, Cut locus, Conjugate locus}
\label{a-conj}

Consider a trivializable $(n,m)$-sub-Riemannian manifold. Solutions to the optimal control problem 
\r{sopra} are computed via the Pontryagin Maximum Principle (PMP for short, see for instance \cite{agra-book,libro,jurd-book,pontlibro}) that is a first order
necessary condition for optimality and generalizes the Weierstra\ss \
conditions of Calculus of Variations. For each optimal curve, the PMP
provides a lift to the cotangent bundle that is a solution to a suitable
pseudo--Hamiltonian system.\\ 

\bt[Pontryagin Maximum Principle for the problem \r{sopra}]
\it Let $M$ be a $n$-dimensional smooth manifold and consider the minimization problem \r{sopra}, in the class of Lipschitz continuous curves, 
 where  $F_i$, $i=1,\ldots, m$ are smooth vector fields on $M$ and the final time $T$ is fixed.  Consider  the map $\mathscr{H}:T^\ast M\times\R\times \R^m\to \R$ defined  by
\bqn
\mathscr{H}(q,\lam,p_0,u)&:=&<\lam,\sum_{i=1}^m u_i F_i(q)>+p_0 \sum_{i=1}^m u_i^2(t).
\eqnn
If the curve $q(.):[0,T]\to M$ corresponding to the control $u(.):[0,T]\to\R^m$ 
is optimal then there exist a {never vanishing}
Lipschitz continuous {covector}
$\lam(.):t\in[0,T]\mapsto \lam(t)\in
T^\ast_{q(t)}M$ and a constant $p_0\leq 0$ such that, for a.e. $t\in 
[0,T]$:
\begin{description}
\item[(i)]
$\dot q(t)=\displaystyle{\frac{\partial \mathscr{H}
}{\partial \lam}(q(t),\lam(t),p_0,u(t))}$,
\item[(ii)] $\dot \lam(t)=-\displaystyle{\frac{\partial \mathscr{H}
}{\partial q}(q(t),\lam(t),p_0,u(t))}$,
\item[(iii)] $\frac{\partial\mathscr{H} }{\partial u}(q(t),\lam(t),p_0,u(t))=0.$
\ed
\et
\brem
A curve $q(.):[0,T]\to M$ satisfying the PMP is said to be an {\it extremal}. In general, an extremal may correspond to more than one pair $(\lam(.),p_0)$. 
If an extremal satisfies the PMP  with $p_0\neq 0 $, then it is called a {\it normal extremals}. If it satisfies the PMP with $p_0= 0 $ it is called an {\it abnormal extremal}. An extremal can be both normal and abnormal. For normal extremals one can normalize $p_0=-1/2$.

If an extremal satisfies the PMP only with $p_0=0$, then it is called a {\it strict abnormal extremal}. If a strict abnormal extremal is optimal, then  it is called {\it a strict abnormal minimizer}. For a deep analysis of abnormal extremals in sub-Riemannian geometry, see \cite{bonnard,trelat}.
\erem

It is well known that all normal extremals are geodesics  (see for instance \cite{agra-book}). Moreover if there are no strict abnormal minimizers then all geodesics are normal extremals for some fixed final time $T$.
This is the case for the so called 3-D {\it contact case}, i.e. a  sub-Riemannian manifold of dimension 3 for which $\Delta=\span{F_1,F_2}$ where $F_1,F_2$ is a pair of vector fields such that for all $q\in M$, $\span{F_1(q),F_2(q), [F_1(q),F_2(q)]}=T_q M$. Indeed for contact structures there are no abnormal extremals (even non strict).

In this case, from {\bf (iii)} one gets $u_i(t)=<\lam(t),F_i(t)>$, $i=1\ldots ,m$ and the PMP becomes  much simpler:   a curve $q(.)$ is a geodesic if and only if it is the projection on $M$ of a 
solution $(\lam(t),q(t))$ for the Hamiltonian system on $T^\ast M$ corresponding to
\bqn
H(\lam,q)=\frac12(\sum_{i=1}^m<\lam,F_i(q)>^2),~~~~
\mbox{$q\in M$, $\lam\in T^\ast_q M$}.
\eqnn
satisfying $H(\lam(0),q(0))\neq 0$.


\brem Notice that $H$ is constant along any given solution of the Hamiltonian  system. Moreover, $H=\frac{1}{2}$ if and only if the geodesic is parameterized by arclength. In the following, for simplicity of notation, we assume that all geodesics are defined for $t\in [0,+\infty)$.
\erem
Fix $q_0\in M$. For every $\lam_0\in T_{q_0}^\ast M$ satisfying
\bqn
H(\lam_0, q_0)=1/2
\eqnl{1/2}
and every $t>0$ define the {\it exponential map} $\Exp(\lam_0,t)$ as the projection on $M$ of the solution, 
evaluated at time $t$, of the Hamiltonian system associated with $H$,  with initial condition $\lam(0)=\lam_0$ and $q(0)=q_0$.
Notice that condition \r{1/2} defines a hypercylinder $\Lambda_{q_0}\simeq S^{m-1}\times\R^{n-m}$ in  $T_{q_0}^\ast M$.
\bdeff
The {\bf conjugate locus from $q_0$} is the set $C_{q_0}$ of critical values of the map \mmfunz{\EXP}{\Lam_{q_0}\times \R^+}{M}{(\lam_0,t)}{\Exp(\lam_0,t).}

For every $\bar\lam_0\in\Lam_{q_0}$, let $t(\bar\lam_0)$ be the n-th positive time, if it exists, for which the map $(\lam_0,t)\mapsto \EXP(\lam_0,t)$ is singular at $(\bar\lam_0,t(\bar\lam_0))$. The {\bf n-th conjugate locus} from $q_0$ $~~~C^n_{q_0}$ is the set $\{\EXP(\bar\lam_0,t(\bar\lam_0))\mid t(\bar\lam_0) \mbox{ 
exists}\}$.

The {\bf cut locus} from $q_0$ is the set $K_{q_0}$ of points reached optimally by more than one geodesic, i.e., the set
\bqn
K_{q_0}=\Pg{q\in M \mid \exists~\lam_1,\lam_2\in\Lambda_{q_0},~\lam_1\neq\lam_2,~t\in\R^+\mbox{ such that~ }\ba{l}  q=\EXP(\lam_1,t)=\EXP(\lam_2,t),\mbox{~~~and~~~}\\ \EXP(\lam_1,\cdot),
\EXP(\lam_2,\cdot)\mbox{ optimal in }[0,t].\ea}
\eqnn
\edeff
\brem
It is a standard fact that for every $\bar\lam_0$ satisfying \r{1/2}, the set $T(\bar\lam_0)=\{\bar t>0\mid 
\EXP(\lam,t)$ is singular at $(\bar\lam_0,\bar t) \}$ is a discrete set (see for instance \cite{agra-book}).
\erem
\brem
\label{rem-agra}
Let $(M,\Delta,\g)$ be a sub-Riemannian manifold. Fix $q_0\in M$ and assume: {\bf (i)} each point of $M$ is reached by an optimal geodesic  starting from $q_0$; {\bf (ii)}  there are no abnormal minimizers.
The following facts are well known (a proof in the 3-D contact case can be found in \cite{agra-exp}). 
\bi
\i the first conjugate locus $C^1_{q_0}$ is the set of points where the geodesics starting from $q_0$   lose local optimality;
\i if $q(.)$ is a geodesic starting from $q_0$ and $\bar t$ is the first positive time such that $q(\bar t)\in K_{q_0}\cup C^1_{q_0}$, then $q(.)$ loses optimality in $\bar{t}$, i.e. it is optimal in $[0,\bar{t}]$ and not optimal in $[0,t]$ for any $t>\bar{t}$;
\i if a geodesic $q(.)$ starting from $q_0$ loses optimality at $\bar t>0$, then $q(\bar t)\in K_{q_0}\cup C^1_{q_0}$;
\ei
As a consequence, when the first conjugate locus is included in the cut locus, the cut locus is the set of points where the geodesics lose optimality.
\erem
\brem
It is well known that, while in Riemannian geometry $K_{q_0}$ is never adjacent to $q_0$, in sub-Riemannian geometry this is always the case. See \cite{Agra-sub}.
\erem

\subsection{$\k\oplus\p$ sub-Riemannian manifolds}
\llabel{s-k+p}
Let $\l$ be a simple Lie algebra and $\Kil(X,Y)=Tr(ad_X\circ ad_Y)$ its Killing 
form. Recall that the Killing form 
defines a non-degenerate pseudo scalar product on $\l$.
In the following we recall what
we mean by a Cartan decomposition of $\l$.      
\bdeff 
A Cartan decomposition of a simple Lie algebra
$\l$ is any
decomposition of the form:
\bqn
\l=\k\oplus\p, 
 \mbox{ where } [\k,\k]\subseteq\k,~~ [\p,\p]\subseteq\k,~~
[\k,\p]\subseteq\p.  
\eqnn 
\edeff
\bdeff 
Let $G$ be a simple Lie group with Lie algebra $\l$. Let $\l=\k\oplus \p$ be a Cartan decomposition of $\l$. In the case in which $G$ is noncompact assume that 
$\k$ is the maximal compact subalgebra of $\l$.

On $G$, consider the distribution  $\Delta(g)=g\p$ endowed with the Riemannian metric ${\mathbf g}_g(v_1,v_2)=<g^{-1}v_1,g^{-1}v_2>$ where $<~,~>:=\al~\Kil\big|_\p(~,~)$ and $\al<0$ (resp. $\al>0$) if $G$ is compact (resp. non 
compact).

In this case we say that $(G, \Delta, {\mathbf g})$ is a $\k\oplus \p$ sub-Riemannian manifold.
\llabel{d-k+p-problem}
\edeff 

The constant $\al$ is clearly not relevant. It is chosen just to obtain
good normalizations.
\brem
In the compact (resp. noncompact) case the fact that ${\mathbf g}$ is positive definite on 
$\Delta$ is guaranteed  by the requirement $\al<0$ (resp. by the requirements $\al>0$  and $\k$ 
maximal compact 
subalgebra).
\erem

Let $\{X_j\}$ be an orthonormal frame for the subspace $\p\subset\l$, with respect to the
metric defined in Definition \ref{d-k+p-problem}. 
Then the problem of finding the minimal energy between the identity and a point $g_1\in G$ in fixed time $T$ becomes the left-invariant optimal control problem
\bqn
\dot g=g\left(\sum_j u_jX_j\right),~~~~u_j\in L^\infty(0,T)\,,
\int_{0}^{T}\sum_ju_j^2(t)~dt\to\min,~~g(0)=\Id,~~~g(T)=g_1.
\eqnn
\brem \llabel{rem-minesiste}
This problem admits a solution, see for instance Chapter 5 of \cite{piccoli}.
\erem

For $\k\oplus\p$ sub-Riemannian manifolds, one can prove that strict abnormal extremals are never optimal, since the {\it Goh condition} (see \cite{agra-book}) is never satisfied. Moreover, the Hamiltonian system given by the Pontryagin Maximum Principle is integrable and the explicit expression of geodesics starting from $g_0$ and parameterized by arclength is
\bqn
g(t)=g_0 e^{(A_k+A_p)t}e^{-A_k t},
\eqnn
where $A_k\in\k$, $A_p\in\p$ and $<A_p,A_p>=1$. This formula is known from long time in the community. It was used  independently by Agrachev \cite{agra-ICM}, Brockett \cite{brokko-cdc99} and Kupka (oral communication). The first complete proof was written by Jurdjevic in \cite{jurd-MCT}. The proof that strict  abnormal extremals are never optimal was first written  in \cite{q2}.  See also  \cite{agra-book,jurd-SCL}.\\\\

{\bf Acknowledgements:} We thank A. Agrachev for his help in recognizing the structure of the projective tangent bundle. We thank L. Paoluzzi for many explanations on lens spaces. We thank M. Saponi and G. Senaldi for their help with pictures.

\end{document}